\documentclass[12 pt]{article}
\usepackage[all]{xy}
\usepackage{pdflscape}
\usepackage[hypcap=false]{caption}
\usepackage{amsmath,amssymb,amsfonts,latexsym,float,graphics,epsfig} 
\usepackage{setspace}
\usepackage{xcolor}
\usepackage[top=2cm, bottom=2cm, left=2cm, right=2cm]{geometry}
\usepackage[colorlinks]{hyperref}
\usepackage{tikz}
\usetikzlibrary{cd}  
\newtheorem{theorem}{Theorem}[section]
\newtheorem{corollary}[theorem]{Corollary}

\title{On Product Lie Algebroids, and Collective Motion}

\begin{document}

\maketitle

\begin{center} 
Begüm Ateşli$^{\ast,\ast\ast,}$\footnote{e-mails: 
\href{mailto:b.atesli@gtu.edu.tr}{b.atesli@gtu.edu.tr}, \href{mailto:begumatesli@itu.edu.tr}{begumatesli@itu.edu.tr} corresponding author,}, Oğul Esen $^{\ast,\dagger,}$\footnote{e-mail: 
\href{oesen@gtu.edu.tr}{oesen@gtu.edu.tr},} and Serkan Sütlü$^{\ast,}$\footnote{e-mail: 
\href{serkansutlu@gtu.edu.tr}{serkansutlu@gtu.edu.tr},}\\

\bigskip
$^\ast$Department of Mathematics, \\ Gebze Technical University, 41400 Gebze,
Kocaeli, Turkey.
\bigskip

$^{\ast\ast}$Department of Mathematics Engineering \\ İstanbul Technical University,  34467 Maslak, İstanbul

\bigskip

$^\dagger$Center for Mathematics and its Applications,\\  Khazar University, Baku, AZ1096,
Azerbaijan

\begin{abstract}
This work explores the geometrical/algebraic framework of Lie algebroids, with a specific focus on the decoupling and coupling phenomena within the bicocycle double cross product realization. The bicocycle double cross product theory serves as the most general method for (de)coupling an algebroid into the direct sum of two vector bundles in the presence of mutual \textit{representations}, along with two twisted cocycle terms. Consequently, it encompasses unified product, double cross product (matched pairs), semi-direct product, and cocycle extension frameworks as particular instances. In addition to algebraic constructions, the research extends to both reversible and irreversible Lagrangian and Hamiltonian dynamics on (de)coupled Lie algebroids, as well as Euler-Poincar\'{e}-(Herglotz) and Lie-Poisson-(Herglotz) dynamics on (de)coupled Lie algebras, providing insights into potential physical applications.

\smallskip

\smallskip


\end{abstract}

\end{center}

\tableofcontents
\onehalfspacing

\setlength{\parindent}{2em}
\setlength{\parskip}{3ex}

\section{Introduction}

The problem addressed in the present manuscript, which may be referred to as the ``(de)coupling problem'', is to obtain the equations of motion that determine the collective behavior of two systems in a mutual action-reaction relationship. In case two dynamical systems or two characters of a physical phenomenon (like Magnetohydrodynamics, where one considers both Maxwellian and fluid character) are coupled, they intervene in each other's independent behavior. As such, the equations of coupled (matched) motion of these two systems/characters cannot be achieved by simply juxtaposing the equations of motions of the individual systems/characters. On the contrary, due to the mutual interaction, the equations that govern the collective motion are expected to contain terms that cannot be retrieved by the equations of the individual systems/characters. 

The (de)coupling process can be examined in two steps. One first needs to determine the extensions of the configuration spaces in a pure algebraic/geometric way. Then, the equations of motion can be coupled by the extensions obtained in the first step. We shall hereby confine ourselves to Lie algebroids and the associated dynamical equations within this geometric framework \cite{MackenzieDG,Mackenzie-book}.   
Our goal, accordingly, is to establish Lie algebroid extensions from a purely algebraic perspective while at the same time maintaining a close relationship with physical applications from both Lagrangian and Hamiltonian viewpoints.

There are several works in the present literature that have addressed the (de)coupling problem in different levels. 
A quick inspection of the algebraic/geometric foundations of certain physical systems (such as Maxwell-Vlasov dynamics, magnetohydrodynamics, and compressible fluid dynamics) reveals that the underlying configuration spaces are \emph{semi-direct product Lie groups} \cite{CeHoMaRa98,kuper83,marsden83b,marsden1984semidirect}. In such physical systems, the reduced Lagrangian (Euler-Poincar\'e) dynamics is available on the \emph{semi-direct product Lie algebra}, whereas the reduced Hamiltonian (Lie-Poisson) dynamics are on the dual space. 
The semi-direct product theory is an example of a Lie group/algebra extension theory that allows only a one-sided (not mutual) action. 

There are generalizations of the semi-direct product theory.  \emph{Double cross product} (matched pair) theory is one such extension, wherein the mutual interactions (actions) of two Lie groups/algebras are allowed \cite{Ma90,majid1990physicsa}. This theory has found applications in physics either. Among examples are the matched pair Lagrangian and Hamiltonian dynamics studied in \cite{esen2016hamiltonian,esen2017lagrangian}, the higher order Lagrangian dynamics in \cite{EsenKudeSutlu21}, and the dissipative systems in \cite{esen2021extensions} and the discrete dynamics in \cite{esen2018matched}. 
These results then soon applied to areas such as electromagnetic theory \cite{Esvd17}, kinetic moments of plasma dynamics, and the algebraic relationship between plasma and fluid \cite{EsGrMiGu19,EsSu21}, as well as the moments in chemical kinetics \cite{AjChEsGrKlPa}. 

We shall further refer the interested reader to \cite{Mokr97} for a detailed discussion on double cross products (matched pairs) within the Lie algebroid framework. 

A more general extension theory (than double cross products) is known as the \emph{unified product theory} \cite{AgorMili11,Agore-Lie,AgorMili14-II,AgoreMilitaru-book}. From the point of view of coupling, and in the level of Lie algebra(oids), unified products correspond to compose a Lie algebra(oid) with a vector bundle to construct a Lie algebra(oid). This is achieved in the presence of a \emph{twisted cocycle}. As such, a unified product has also been called a \emph{cocycle double cross product} in \cite{esen2021bicocycle}. 
From the decomposition perspective, a unified product (a cocycle double cross product) may be interpreted as a recipe to study the quotient space of an algebroid by a subalgebroid of it. 
In the Lie algebra setting, we cite a recent work \cite{uccgun2024dynamics} where unified products are used to study the dynamics on homogeneous spaces.  

Unified products (cocycle double cross products), subsume two independent avenues of extensions; the double cross products, and 2-cocycle extensions. They thus offer much more flexibility in (de)coupling a geometric structure, or a physical phenomena, compared to what might be offered by double cross products or 2-cocycle extensions.

Despite this pliability, the absence of applications (of unified products to Lagrangian and Hamiltonian dynamics as well as any interconnection with physical phenomena) in the literature stands as a strong motivation for us. 

The cocycle double cross products have been upgraded recently in \cite{esen2021bicocycle} to a level that requires the use of two twisted cocycles. The new class of extension was accordingly called \emph{bicocycle double cross product}, which shall occasionally be abbreviated as BDCP. In the level of Lie algebroids, a BDCP corresponds to a couple of two vector bundles arriving at a Lie algebroid on the direct sum vector bundle. The non-closure of the (induced) brackets on constitutive vector bundles is accommodated by two (twisted) cocycles, as a manifestation of the name. 

BDCPs, by their very nature, thus render the ultimate flexibility for (de)coupling. However, the literature is void in applications, providing another motivation for the present manuscript.

We shall continue with a summary of the families of extensions mentioned above (explicit definitions will be given in the main body of the text). We shall next spell out our objectives.

\subsubsection*{The Hierarchy of Extensions}

Let $\mathcal{M}$ be a Lie algebroid over a base manifold $M$, along with two complementary vector subbundles $\mathcal{A}$ and $\mathcal{B}$ so that $\mathcal{M}$ is the Whitney sum\footnote{By a slight abuse of notation, we shall use $\times$ for the Whitney sum of vector bundles over the same base.}, that is,
\[
\mathcal{M} \cong \mathcal{A} \times \mathcal{B}
\]
as vector bundles. The bracket operation $[\bullet,\bullet]_{\mathcal{M}}:\Gamma{(\mathcal{M})} \times \Gamma{(\mathcal{M})} \mapsto \Gamma{(\mathcal{M})}$ on the space of sections of $\mathcal{M}$, then, induces corresponding operations 
\begin{equation}\label{the-psi-map}
   \Gamma{(\mathcal{A})} \times \Gamma{(\mathcal{A})} \to \Gamma{(\mathcal{M})}
\end{equation}
on $\mathcal{A}$ and  
\begin{equation}\label{the-zeta-map}
   \Gamma{(\mathcal{B})} \times \Gamma{(\mathcal{B})} \to \Gamma{(\mathcal{M})}
\end{equation}
on $\mathcal{B}$. It follows at once that neither \eqref{the-psi-map} nor \eqref{the-zeta-map} are a priori closed (neither they satisfy readily the Jacobi or Leibniz identities). Failure in closure of the induced operations are compensated by two mappings
\begin{equation}\label{the-psi-map-2}
   \psi: \Gamma{(\mathcal{B})} \times \Gamma{(\mathcal{B})} \to \Gamma{(\mathcal{A})}
\end{equation}
and  
\begin{equation}\label{the-zeta-map-2}
  \zeta: \Gamma{(\mathcal{A})} \times \Gamma{(\mathcal{A})} \to \Gamma{(\mathcal{B})},
\end{equation}
which are usually named (twisted) cocycles, in addition to
\begin{equation}\label{the-phi-map}
   \phi: \Gamma{(\mathcal{A})} \times \Gamma{(\mathcal{A})} \to \Gamma{(\mathcal{A})}
\end{equation}
on $\mathcal{A}$ and  
\begin{equation}\label{the-theta-map}
  \theta: \Gamma(\mathcal{B}) \times \Gamma(\mathcal{B}) \to \Gamma(\mathcal{B})
\end{equation}
on $\mathcal{B}$.
 
On the other hand, it follows from the direct (Whitney) sum decomposition of $\mathcal{M}$, the mixed bracket (within $\mathcal{M}$, of an element in $\mathcal{B}$ and an element in $\mathcal{A}$) yields mutual \emph{representations}\footnote{To be more precise, two mappings $\Gamma{(\mathcal{B})} \times \Gamma{(\mathcal{A})} \to \Gamma{(\mathcal{A})}$ and $\Gamma{(\mathcal{B})} \times \Gamma{(\mathcal{A})} \to \Gamma{(\mathcal{B})}$, which are called \emph{weak actions}. To represent the weak actions, we shall make use of the italic versions of the notations we adopt for the regular actions.} of $\mathcal{A}$ and $\mathcal{B}$ on each other. 

The bicocycle double cross product Lie algebroid $\mathcal{M} \cong \mathcal{A} \times \mathcal{B}$ is then denoted by $\mathcal{A} {\hspace{.1cm}}_{\zeta\hspace{-.1cm}} \bowtie _\psi\mathcal{B}$, where the indices $\psi$ and $\zeta$ highlight the (twisted) cocycles \eqref{the-psi-map-2} and \eqref{the-zeta-map-2}, and the triangles $\triangleright$ and $\triangleleft$ emphasize the mutual (weak) actions. While we will study this more explicitly in the main body of the paper, to clarify the hierarchical structure that follows, let us write the bicocycle double cross product Lie bracket:
\begin{equation}
	\begin{split}
 	\Gamma(\mathcal{A}  \times  \mathcal{B}) \times  	\Gamma(\mathcal{A}  \times \mathcal{B}) & \longrightarrow 	\Gamma(\mathcal{A} \times  \mathcal{B}), \\		
		[(X_1,Y_1),(X_2,Y_2)]_{\mathcal{M}} & = \Big(\phi(X_1,X_2)- \rho(Y_2,X_1) + \rho(Y_1,X_2) + \psi(Y_1,Y_2) ,\\
		&\qquad \qquad \theta(Y_1,Y_2) - \sigma(X_2,Y_1) + \sigma(X_1,Y_2)+\zeta(X_1,X_2)\Big)
	\end{split}
\end{equation}
in terms of the mappings in \eqref{the-psi-map-2}-\eqref{the-theta-map} as well as the weak action $\rho$ of $\mathcal{B}$ on $\mathcal{A}$, and the weak action $\sigma$ of $\mathcal{A}$ on $\mathcal{B}$. 

In the case one of the (twisted) cocycles $\psi$ in \eqref{the-psi-map-2} or $\zeta$ in \eqref{the-zeta-map-2} is trivial (that is identically zero), the Lie algebroid $\mathcal{M}$ becomes a unified product (cocycle double cross product) of $\mathcal{A}$ and $\mathcal{B}$, and is denoted by $\mathcal{A} {\hspace{.1cm}}_{\zeta\hspace{-.1cm}} \bowtie \mathcal{B}$ if \eqref{the-psi-map-2} is trivial, and $\mathcal{A} \bowtie _\psi\mathcal{B}$ if \eqref{the-zeta-map-2} is trivial. 

Let us now consider the unified product (cocycle double cross product) $\mathcal{M}\cong\mathcal{A} \bowtie _\psi\mathcal{B}$. In this case the weak action $\sigma$ of $\mathcal{A}$ on $\mathcal{B}$ becomes an action.

If, further, the (twisted) cocycle $\psi$ in \eqref{the-psi-map-2} is trivial either, then $\mathcal{M}$ happens to be a double cross product of $\mathcal{A}$ and $\mathcal{B}$, and is denoted by $\mathcal{M}\cong\mathcal{A} \bowtie \mathcal{B}$. In this case, the mutual weak actions of $\mathcal{A}$ and $\mathcal{B}$ on each other are both genuine actions, and $\mathcal{A}$ and $\mathcal{B}$ are both Lie subalgebroids of $\mathcal{M}$.

If, on the other extreme, the action of $\mathcal{A}$ on $\mathcal{B}$ is trivial, along with the bracket structure on the space $\Gamma(\mathcal{A})$ of sections of $\mathcal{A}$, then the product Lie algebroid $\mathcal{M}$ turns out to be a 2-cocycle extension of $\mathcal{A}$ by $\mathcal{B}$. In this case we use the notation $\mathcal{M} \cong \mathcal{A} \rtimes_\psi \mathcal{B}$.

It follows from the assumptions of this case that $\mathcal{B}$ happens to be a Lie algebroid, and the weak action of $\mathcal{B}$ on $\mathcal{A}$ is now a proper action. Further, the mapping \eqref{the-psi-map-2} turns out to be a 2-cocycle in the Lie algebroid cohomology of $\mathcal{B}$, with coefficients in the vector bundle $\mathcal{A}$. However, $\mathcal{B}$ is not a Lie subalgebroid of $\mathcal{M}$, while $\mathcal{A}$ is (albeit a trivial one). Let us note also that if, further, the action of $\mathcal{B}$ on $\mathcal{A}$ is assumed to be trivial, then $\mathcal{A}$ becomes central (as a Lie subalgebroid) in $\mathcal{M}$, and the 2-cocycle extension $\mathcal{M} \cong \mathcal{A} \rtimes_\psi \mathcal{B}$ becomes a central extension.

The common ground of the 2-cocycle extensions and double cross products is given by \emph{semi-direct products}. More precisely, if the cocycle \eqref{the-psi-map-2} is assumed to be trivial in the 2-cocycle extension $\mathcal{M} \cong \mathcal{A} \rtimes_\psi \mathcal{B}$, or if the action of $\mathcal{A}$ on $\mathcal{B}$ is taken to be trivial in the double cross product $\mathcal{M} \cong \mathcal{A} \bowtie \mathcal{B}$, then the object to arrive is a semi-direct product built solely on the action of $\mathcal{B}$ on $\mathcal{A}$, which is denoted by $\mathcal{M} \cong \mathcal{A} \rtimes \mathcal{B}$.

The whole hierarchy is illustrated in the following diagram, within which we use the abbreviation DC for ``Double Cross''.

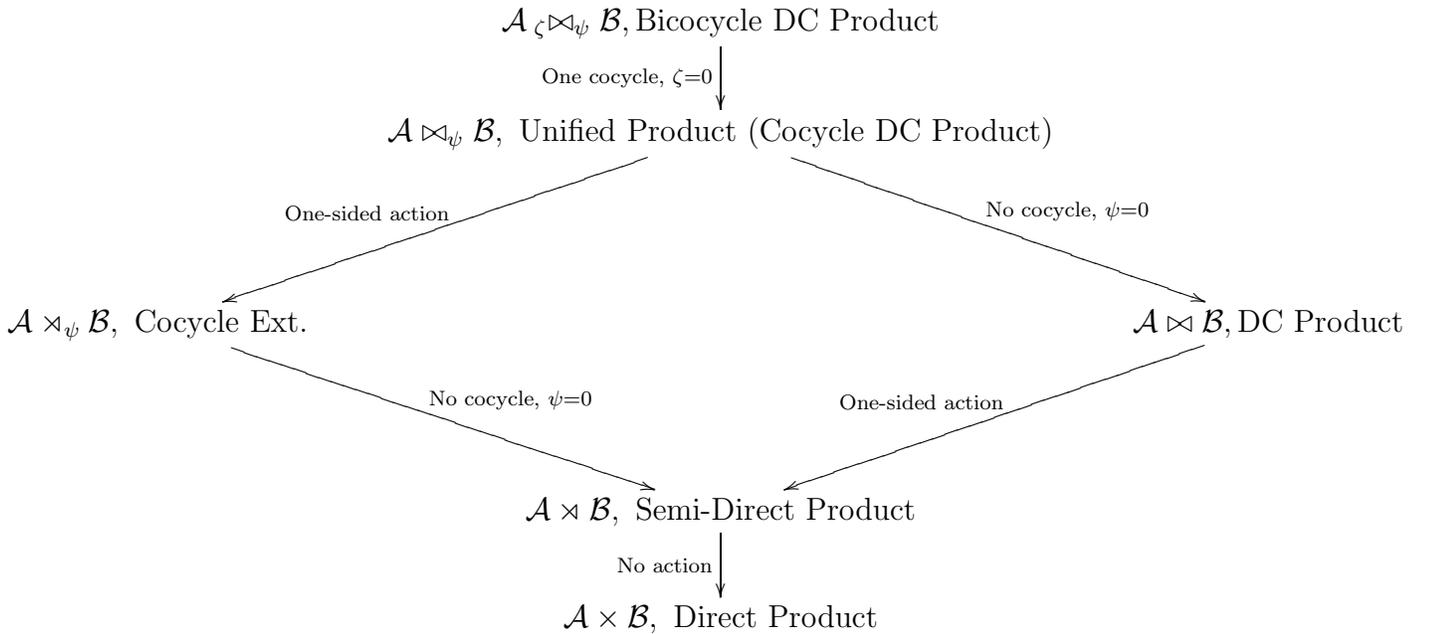
 \begin{equation*}\label{Diagram}
	\xymatrix{ 
 & \mathcal{A} {\hspace{.1cm}}_{\zeta\hspace{-.1cm}} \bowtie _\psi\mathcal{B}, \text{Bicocycle DC Product}\ar[d]_{\text{One cocycle, }\zeta=0}
 \\ &  \mathcal{A}\bowtie _\psi\mathcal{B}, \text{ Unified Product (Cocycle DC Product)} \ar[ddl]_{\text{One-sided  action}\qquad}\ar[ddr]^{\qquad\text{No cocycle, }\psi=0}\\
 & & & \\  \mathcal{A}\rtimes _\psi\mathcal{B},
		\text{ Cocycle Ext.} 
		\ar[ddr]^{\qquad\text{No cocycle, }\psi=0} &     &\mathcal{A}\bowtie\mathcal{B}, \text{DC Product} \ar[ddl]_{\text{One-sided  action}\qquad}\\
  & & & \\
 & \mathcal{A}\rtimes \mathcal{B},
 \text{ Semi-Direct Product} \ar[d]_{\text{No action}}
 \\
  & \mathcal{A}\times \mathcal{B}, \text{ Direct Product} 
 }
\end{equation*}

We shall now review the objectives of the manuscript, along the lines of the product Lie algebroids summarized above.

\subsubsection*{Goal I: Bicocycle Double Cross Product Lie Algebroids.}

Various components of the cross product hierarchy above have been developed for a plethora of algebraic / geometric objects (such as Lie groups, Lie algebras, Hopf algebras, etc.), but not for Lie algebroids. More precisely, the very first objective of the present manuscript is to develop the theory of bicocycle double cross product Lie algebroids, namely the top level in the hierarchy of product Lie algebroids.

We then mean to use the (bicocycle double) cross product theory of Lie algebroids in the investigation of dynamical systems. More precisely, the theory we cultivate will allow us to decompose a system modelled by a Lie algebroid into simpler parts. The individual dynamics that are governed by these simpler components, then, unites to shed light to the dynamics of combined system. 

It will thus be proper to include a brief outline of the dynamics on Lie algebroids.

\subsubsection*{Dynamics on Lie Algebroids.} 

Given a Lie algebroid $(\mathcal{A},\tau,M,\mathfrak{a}_{\mathcal{A}},[\bullet,\bullet])$, the Hamiltonian dynamics exist on the dual bundle $\mathcal{A}^*$, \cite{de2005lagrangian,GrabowskaAlg}. Assume a local coordinate system $(x^i)$ on $M$, along with $(x^i,y^\alpha)$ on $\mathcal{A}$. 
In local (dual) coordinates $(x^i,y_\alpha)$ on $\mathcal{A}^*$, given a Hamiltonian function $\mathcal{H}=\mathcal{H}(x^i,y_\alpha)$, one has \emph{Hamilton's equations}
\begin{equation}\label{ham-eq-int}
 \frac{d x^i}{dt} = (\mathfrak{a}_{\mathcal{A}})_{\alpha }^{i}\frac{\partial \mathcal{H}}{\partial y_{\alpha}},\qquad
 \frac{d y_{\alpha}}{dt} = - \mathfrak{H}_{\alpha \beta}^{\gamma}y_{\gamma} \frac{\partial \mathcal{H}}{\partial y_{\beta}} - (\mathfrak{a}_{\mathcal{A}})_{\alpha }^{i}\frac{\partial \mathcal{H}}{\partial x^i},
\end{equation}
where $(\mathfrak{a}_{\mathcal{A}} )^i_\alpha$ is the local realization of the anchor map while $\mathfrak{H}^\gamma_{\alpha\beta}$ denotes the structure functions of the Lie algebroid. 

If the Lie algebroid is a tangent bundle, \eqref{ham-eq-int} corresponds to the classical Hamilton's equation on the cotangent bundle, arising from the existence of the canonical symplectic structure. 

If, on the other hand, the Lie algebroid is a Lie algebra $\mathfrak{g}$ (the base manifold being a singleton, the anchor map is trivial), \eqref{ham-eq-int} represents the \emph{Lie-Poisson equations} on the dual space $\mathfrak{g}^*$. More precisely, assuming the local coordinates $(y_\alpha)$ on $\mathfrak{g}^*$, and considering a Hamiltonian function $h=h(y_\alpha)$, the Lie-Poisson equations are given by 
\begin{equation}\label{ham-eq-LP-int}
 \frac{d y_{\alpha}}{dt} = - \mathfrak{H}_{\alpha \beta}^{\gamma}y_{\gamma} \frac{\partial h}{\partial y_{\beta}},
\end{equation}
see for instance \cite{AbMa78,holm2009geometric,MaRa99}. Here, $\mathfrak{H}_{\alpha \beta}^{\gamma}$ are now the structure constants of the Lie algebra. The Lie-Poisson equations, describing the dynamics on $\mathfrak{g}^*$, can also be derived through the reduction by symmetry method. 

The Lie algebroid framework encapsulates the Lagrangian dynamics as well. In local coordinates $(x^i, y^\alpha)$ on a Lie algebroid $\mathcal{A}$, the dynamics encoded by a Lagrangian function $L = L(x^i, y^\alpha)$ is described by the \emph{Euler-Lagrange equations}
\begin{equation} \label{AEL-intro}
\frac{dx^i}{dt} = (\mathfrak{a}_{\mathcal{A}} )^i_\alpha y^\alpha
,\qquad
\frac{d}{dt}\frac{\partial L}{\partial y^\alpha} = (\mathfrak{a}_{\mathcal{A}} )^i_\alpha \frac{\partial L}{\partial x^i} +  \mathfrak{H}^\gamma_{\alpha\beta}y^\beta\frac{\partial L}{\partial y^\gamma},
\end{equation}
see for instance \cite{WeinLag}.

The coordinate-independent realization of the above formulation was later established in \cite{martinez2001}. In the particular case where the Lie algebroid is a tangent bundle, the equations \eqref{AEL-intro} correspond to the classical Euler-Lagrange equations. If, on the other extreme, the Lie algebroid is a Lie algebra $\mathfrak{g}$ with coordinates $(y^\alpha)$, \eqref{AEL-intro} transforms into the \emph{Euler-Poincar\'{e} equations}
\begin{equation} \label{EP-intro} 
\frac{d}{dt}\frac{\partial l}{\partial y^\alpha}=  \mathfrak{H}^\gamma_{\alpha\beta}y^\beta\frac{\partial l}{\partial y^\gamma},
\end{equation}
for a Lagrangian function $l=l(y^\alpha)$, \cite{MaRa99}. 
In regular cases, the connection between Lagrangian and Hamiltonian dynamics is established through the Legendre transformation.

The reversible nature of both \eqref{ham-eq-int} and \eqref{AEL-intro} highlights the limitations in their applicability to irreversible or/and dissipative processes. The literature points out that the contact structures provide a suitable avenue for  irreversible dynamics, \cite{Br17,BrCrTa17,cannarsa2019herglotz,LeLa19,de2020review,LeSa17,esen2021contact,esen2021implicit}. This approach has recently amplified the attention to the field of thermodynamics, \cite{Bravetti19,OMM-1,OMM-2,goshTermo,Goto,Grmela-Contact,Mrugala,Simoes2020b}. 

A generic way to arrive at a contact manifold is to extend a symplectic manifold, say $T^*Q$, by a copy of the real line. In Darboux' coordinates $(q^i,p_i,z)$ on the product bundle $T^*Q\times \mathbb{R}$, for a Hamiltonian function $\mathcal{H}=\mathcal{H}(x^i,p_i,z)$, the \emph{contact Hamilton's equations} are 
\begin{equation}\label{conham-X-intro}
\frac{dx^i}{dt} = \frac{\partial \mathcal{H}}{\partial p_{i}}, 
\qquad
\frac{dp_i}{dt}= -\frac{\partial \mathcal{H}}{\partial x^{i}}- 
p_{i}\frac{\partial \mathcal{H}}{\partial z}, 
\qquad \frac{dz}{dt}= p_{i}\frac{\partial \mathcal{H}}{\partial p_{i}} - \mathcal{H}.
\end{equation}
Accordingly, a suitable geometric framework for contact Lagrangian dynamics (that allows dissipation) may be established through the extension $TQ \times \mathbb{R}$ of the tangent bundle $TQ$ by the real line. In induced coordinates $(x^i, \dot{x}^i, z)$ on $TQ \times \mathbb{R}$, $L = L(x^i, \dot{x}^i, z)$ being the Lagrangian function, the corresponding equations
\begin{equation}\label{clagrangian4-intro}
\frac{d}{dt} \frac{\partial L}{\partial \dot{x}^i} - \frac{\partial L}{\partial x^i} =
\frac{\partial L}{\partial \dot{x}^i} \frac{\partial L}{\partial z}, \qquad \frac{dz}{dt} = L(x^i, \dot{x}^i, z)
\end{equation}
are known as \emph{Herglotz equations}, or \emph{generalized Euler-Lagrange equations}, \cite{de2020review}.

The strategy of extending by the real line has recently been applied to Lie algebroids (and their duals), in an attempt to extend the scope of the dynamics on Lie algebroids to a theory that allows dissipation, \cite{anahory2023reduction,anahory2024euler,Simoes2020b}. The scheme is named as the \emph{contactization} (of Lie algebroids).

More explicitly, letting $\mathcal{A}$ be a Lie algebroid, and considering the extension $\mathcal{A}^* \times \mathbb{R}$ of the dual bundle with local coordinates $(x^i, y_\alpha, z)$, the \emph{dissipative Hamilton's equations}, generated by a Hamiltonian function $\mathcal{H} = \mathcal{H}(x^i, y_\alpha, z)$, take the form
\begin{equation}\label{conalghameq-intro} 
	\frac{dx^i}{dt}  =(\mathfrak{a}_{\mathcal{A}})_{\alpha}^{i}\frac{\partial \mathcal{H}}{\partial y_{\alpha}} ,  \quad \frac{dy_\alpha}{dt} =-  (\mathfrak{a}_{\mathcal{A}})_{\alpha}^{i}\frac{\partial \mathcal{H}}{\partial x^{i}}- \mathfrak{H}_{\alpha \beta}^{\gamma} y_{\gamma}\frac{\partial \mathcal{H}}{\partial y_{\beta}}- y_{\alpha}\frac{\partial \mathcal{H}}{\partial z} ,  \quad  \frac{dz}{dt}=  y_{\alpha}\frac{\partial \mathcal{H}}{\partial y_{\alpha}}-\mathcal{H}   .
\end{equation}
It is worth to note that the above set of equations seamlessly combines the structure of the Hamilton's equations \eqref{ham-eq-int} on the dual bundle $\mathcal{A}^*$ with the contact Hamilton's equations presented in \eqref{conham-X-intro}. Indeed, if in particular the Lie algebroid $\mathcal{A}$ happens to be the tangent bundle $TQ$, then the dissipative Hamilton's equations \eqref{conalghameq-intro} naturally yield the contact Hamiltonian dynamics as expressed in \eqref{conham-X-intro}.

If, at the other extreme, the Lie algebroid $\mathcal{A}$ is a Lie algebra $\mathfrak{g}$, the resulting equations are called the \emph{Lie-Poisson-Herglotz equations} on the extended dual space $\mathfrak{g}^*\times \mathbb{R}$. In coordinates $(y_\alpha, z)$ on $\mathfrak{g}^*\times \mathbb{R}$, with a Hamiltonian function $h = h(y_\alpha, z)$, the Lie-Poisson-Herglotz equations take the form
\begin{equation}\label{conalghameq-intro-ref} 
	 \frac{dy_\alpha}{dt} = - \mathfrak{H}_{\alpha \beta}^{\gamma} y_{\gamma}\frac{\partial h}{\partial y_{\beta}}- y_{\alpha}\frac{\partial h}{\partial z} ,  \quad  \frac{dz}{dt}=  y_{\alpha}\frac{\partial h}{\partial y_{\alpha}}-h   .
\end{equation}

From the Lagrangian point of view, let $L = L(x^i, y^\alpha, z)$ denotes a Lagrangian function on the extended Lie algebroid $\mathcal{A}\times \mathbb{R}$ with local coordinates $(x^i, y^\alpha, z)$. Then the \emph{dissipative Euler-Lagrange equations} encoded by this Lagrangian are
\begin{equation}\label{HerlAlg-intro}
   \frac{d}{dt}\frac{\partial L}{\partial y^{\alpha}}-(\mathfrak{a}_{\mathcal{A}})_{\alpha}^{i}\frac{\partial L}{\partial x^{i}}=\mathfrak{H}_{\alpha \beta}^{\gamma}y^{\beta}\frac{\partial L}{\partial y^{\gamma}} +\frac{\partial L}{\partial z}\frac{\partial L}{\partial y^{\alpha}},\quad \frac{dx^i}{dt} = (\mathfrak{a}_{\mathcal{A}})_{\alpha}^{i}y^{\alpha} ,  \quad \frac{dz}{dt}=L(x^i,y^\alpha,z). 
\end{equation} 
Let us note once again that if, on the one extreme, the Lie algebroid $\mathcal{A}$ is the tangent bundle $TQ$, then the dissipative Euler-Lagrange equations \eqref{HerlAlg-intro} turn into Herglotz equations \eqref{clagrangian4-intro} of contact Lagrangian dynamics on the extended tangent bundle $TQ \times \mathbb{R}$. If, on the other extreme, the Lie algebroid $\mathcal{A}$ is in particular a Lie algebra $\mathfrak{g}$, then \eqref{clagrangian4-intro} transform into the \emph{Euler-Poincar\'{e}-Herglotz equations} the Euler-Poincar\'{e}-Herglotz equations:
\begin{equation}\label{HerlAlg-intro-EP}
   \frac{d}{dt}\frac{\partial l}{\partial y^{\alpha}}
   =\mathfrak{H}_{\alpha \beta}^{\gamma}y^{\beta}\frac{\partial l}{\partial y^{\gamma}} +\frac{\partial l}{\partial z}\frac{\partial l}{\partial y^{\alpha}} ,  \quad \frac{dz}{dt}=l(x^i,y^\alpha,z)
\end{equation} 
of a Lagrangian function $l = l(y^\alpha, z)$ on $\mathfrak{g}\times \mathbb{R}$ with coordinates $(y^\alpha, z)$.

In the context of the present work, the formulations \eqref{conalghameq-intro} and \eqref{HerlAlg-intro} hold particular importance due to their extension / product nature, in the level of Lie algebroids. It will hence be well-suited for our second objective; namely to apply the Lie algebroid cross product theory we develop to dynamics on Lie algebroids.

\subsubsection*{Goal II: Application. Bicocycle Double Cross Product Dynamics.}

The second objective of this work is to present a comprehensive framework for both Hamiltonian and Lagrangian dynamics in view of the bicocycle double cross product Lie algebroids outlined above. More explicitly, we shall interpret a Lie algebroid (for the Lagrangian dynamics) and its dual (for the Hamiltonian dynamics) as a bicocycle double cross product in order to recast the equations of motion in terms of the individual equations of motions of the (simpler) pieces of the decomposition. The new formulation of the equations of motion involve possibly some terms that do not arise from these individual equations of motions. They, instead, captures the mutual interaction of the pieces which are subject to the decomposition. Such a reformulation of the equations of motion will indeed involve more number of terms than the original equations. However, the terms that arise from the cross product structure will, in practice, be much easier to deal with.

As far as the Hamiltonian mechanics is concerned, this recipe will allow us to decompose the Hamilton's equations \eqref{ham-eq-int}, and their dissipative counterparts \eqref{conalghameq-intro}. For the Lagrangian mechanics, on the other hand, we shall thus obtain a decomposition of the Euler-Lagrange equations \eqref{AEL-intro}, as well as the dissipative Euler-Lagrange equations \eqref{HerlAlg-intro} for the dissipative scenario. Needless to say, particular instances of these equations (such as the Lie-Poisson-Herglotz equations and Euler-Poincar\'{e}-Herglotz equations) will be regenerated at once. 

Accordingly, this work may be considered as a sequel to \cite{Atesli22,esen2021bicocycle,uccgun2024dynamics} where the bicocycle double cross products were treated in the level of Lie algebras.

\subsubsection*{Outline.} 

This text consists of two main parts addressing the two main objectives presented above. 

Section \ref{sec-unfLiealg} represents the algebraic foundations of the diagram \eqref{Diagram} of the hierarchy of product Lie algebroids. It is in this section that we shall accomplish to obtain the compatibility conditions on two vector bundles whose product (Whitney sum) is a Lie algebroid. The main theorems on bicocycle double cross products will be collected in Subsection \ref{sec-bicocycle}. We next provide two examples. First, in Subsection \ref{Dirac-sec}, we study the decomposition of twisted Poisson geometry in the Dirac bundle setting. Then, in Subsection \ref{decomsec}, we turn our attention to decompose a Jacobi Lie algebroid into a Poisson Lie algebroid and a line bundle. This has been achieved in the level of Lie algebras in \cite{Agore-Jacobi}.

We present the equations of motion, both in Lagrangian and Hamiltonian frameworks and both for reversible and irreversible dynamics, in Section \ref{dyndyn}. More explicitly, in Subsection \ref{dyndyn-1} we recall the equations of motion for the reversible Lagrangian dynamics on a Lie algebroid, as well as those for the reversible Hamiltonian dynamics on the dual bundle. We then regenerate these equations in Subsection \ref{sec-dyn-bicocycle} in the presence of a bicocycle double cross product decomposition of a Lie algebroid (for the Langangian perspective) and its dual (for the Hamiltonian perspective). Subsection \ref{cont-sec} is reserved for a brief summary of the contact dynamics as a dissipative generalization of the classical one in view of the real line extension. Similarly, Subsection \ref{sec-irr-Lie} is withheld to record a brief summary of the dissipative dynamics in the level of Lie algebroids. Finally, in Section \ref{noveldyn} we reproduce the equations for dissipative dynamics in both the Hamiltonian and the Lagrangian settings in view of a bicocycle double cross product decomposition of the Lie algebroid on which the dissipative equations were given.

 \section{Product Lie Algebroids}\label{sec-unfLiealg}

\subsection{Unified Product Lie Algebroids}\label{UPLA-sec}

For the sake of the completeness of the exposition, and for future reference, let us first recall the terminology of Lie algebroids, and then that of the unified product construction in the Lie algebroid framework.

\subsubsection*{Lie Algebroids.} 

Given a manifold $M$, a \emph{Lie algebroid} $\mathcal{A}$ over the base $M$ is a vector bundle $\tau:\mathcal{A} \mapsto M$, equipped with a map 
	\begin{equation} 
\mathfrak{a}_\mathcal{A}:\mathcal{A} \longrightarrow TM
	\end{equation}
 of vector bundles, called the \emph{anchor map}, and a skew-symmetric bilinear bracket $[\bullet,\bullet]_\mathcal{A}$ on the space $\Gamma(\mathcal{A})$ of sections of $\mathcal{A}$ so that
\begin{equation} 
		\begin{split}
				(i)\quad	& \mathfrak{a}_{\mathcal{A}}([X_1,X_2]_\mathcal{A}) = [\mathfrak{a}_{\mathcal{A}}(X_1),\mathfrak{a}_{\mathcal{A}}(X_2)] \\
				(ii)\quad & [X_1,fX_2]_\mathcal{A} = f[X_1,X_2]_\mathcal{A}+\mathcal{L}_{\mathfrak{a}_\mathcal{A}(X_1)}(f)X_2\\
    (iii)\quad & \circlearrowright [X_1, [X_2, X_3]_\mathcal{A}]_\mathcal{A} = 0,
		\end{split}
	\end{equation}
	for any $X_1,X_2,X_3 \in \Gamma(\mathcal{A})$, and any $f\in \mathcal{C}^\infty(M)$ \cite{Mackenzie-book,Mokr97,Para67}. The second equality is called the \emph{Leibniz identity}, while the last one ($\circlearrowright$ refers to the cyclic sum of the relevant elements) is named as the \emph{Jacobi identity}. Accordingly, a Lie algebroid is denoted by a quintuple 
\[
(\mathcal{A},\tau,M,\mathfrak{a}_\mathcal{A},[\bullet,\bullet]_\mathcal{A}),
\]
or in short by $\mathcal{A}$, which may be presented by a commutative diagram 
\[
\xymatrix{
\mathcal{A} \ar[rr]^{\mathfrak{a}_\mathcal{A}} \ar[dr]_\tau &  & TM \ar[dl]^{\tau_M} \\
& M
}
\]
where $\tau_M$ stands for the tangent bundle projection. 

In finite dimensions, let $(x^i)$ be a local coordinate system on $M$. A local coordinate system on $\mathcal{A}$ may hence be denoted by $(x^i, y^\alpha)$. Letting, further, $(e_\alpha)$ be a (projective) basis for the space $\Gamma(\mathcal{A})$ of sections, we obtain the matrix representations of the anchor map $\mathfrak{a}_{\mathcal{A}}$, and the structure functions for the Lie algebroid $\mathcal{A}$ as 
\begin{equation} \label{loc-algebroid}
\mathfrak{a}_{\mathcal{A}} (e_\alpha) = (\mathfrak{a}_{\mathcal{A}})^i_\alpha\frac{\partial}{\partial x^i}, \qquad [e_\alpha,e_\beta]_\mathcal{A} =\mathfrak{H}_{\alpha\beta}^\gamma e_\gamma,
\end{equation}
respectively.

\subsubsection*{Representations of Lie Algebroids.}

The representation of a Lie algebroid $(\mathcal{A},\tau,M,\mathfrak{a}_{\mathcal{A}},[\bullet,\bullet])$ on a vector bundle $(\mathcal{B},\kappa,M)$ is defined to be a map 
\[
\sigma:\Gamma(\mathcal{A})\times \Gamma(\mathcal{B}) \longrightarrow\Gamma(\mathcal{B}), \qquad (X,Y)\mapsto \sigma(X,Y)
\]
subject to
\begin{align*} 
\sigma(fX,Y) &= f\sigma(X,Y) ,\\
\sigma(X,fY) & = f\sigma(X,Y) + \mathcal{L}_{\mathfrak{a}_{\mathcal{A}}(X)}(f)Y,\\
\sigma([X_1,X_2 ]_\mathcal{A},Y) & = \sigma(X_1,\sigma(X_2,Y)) - \sigma(X_2,\sigma(X_1,Y)) 	
\end{align*}
for any $X,X_1,X_2\in\Gamma(A)$, any $Y\in \Gamma(\mathcal{B})$, and any $f \in \mathcal{C}^\infty(M)$.

\subsubsection*{Unified Products of Lie Algebroids.}

Let $(\mathcal{A},\tau,M,\mathfrak{a}_{\mathcal{A}},[\bullet,\bullet]_{\mathcal{A}})$ be a Lie algebroid, and let $(\mathcal{B},\kappa,M)$ be an arbitrary vector bundle together with a vector bundle map $\mathfrak{a}_{\mathcal{B}}:\mathcal{B} \mapsto TM$. 

In an attempt to define a Lie algebroid structure on the Whitney (direct) sum vector bundle $\mathcal{A}\times \mathcal{B}$ over the base manifold $M$, equipped with the anchor
\begin{equation}\label{bobobo}
\mathfrak{a}_{\bowtie_\psi}:=\mathcal{A}\times \mathcal{B} \longrightarrow TM,\qquad 	 (X,Y) \mapsto \mathfrak{a}_{\mathcal{A}}(X) + \mathfrak{a}_{\mathcal{B}}(Y),
\end{equation}
let us assume the existence of two sets of skew-symmetric bilinear maps 
\begin{equation}\label{muratcan}
\psi:\Gamma(\mathcal{B})\times \Gamma(\mathcal{B}) \longrightarrow \Gamma(\mathcal{A}), \qquad \theta:\Gamma(\mathcal{B})\times \Gamma(\mathcal{B}) \longrightarrow \Gamma(\mathcal{B})
\end{equation}
and 
\begin{equation}\label{rolar}
\rho: \Gamma(\mathcal{B})\times \Gamma(\mathcal{A}) \longrightarrow \Gamma(\mathcal{A}), \qquad \sigma: \Gamma(\mathcal{A})\times \Gamma(\mathcal{B}) \longrightarrow \Gamma(\mathcal{B}).
\end{equation}
Then, a natural bracket operation on the space of sections of the product bundle $\mathcal{A}\times \mathcal{B}$ may be given by 
\begin{equation}\label{hohoho}
\begin{split}
\Gamma(\mathcal{A}  \times  \mathcal{B}) \times  	\Gamma(\mathcal{A}  \times \mathcal{B}) & \longrightarrow 	\Gamma(\mathcal{A} \times  \mathcal{B}), \\		
[(X_1,Y_1),(X_2,Y_2)]_{\bowtie_\psi} & := \Big([X_1,X_2]_ \mathcal{A}- \rho(Y_2,X_1) + \rho(Y_1,X_2) + \psi(Y_1,Y_2) ,\\
&\qquad \qquad \theta(Y_1,Y_2) - \sigma(X_2,Y_1) + \sigma(X_1,Y_2)\Big).
\end{split}
\end{equation}
The bracket \eqref{hohoho} is automatically skew-symmetric (as a result of the skew-symmetry of the mappings). As for the Jacobi and Leibniz identities, on the other hand, the mappings of \eqref{muratcan} and \eqref{rolar} are subject to a number of conditions. 
In fact, assuming for the sake of simplicity the latter mapping of \eqref{rolar} is a genuine Lie algebroid action\footnote{This would otherwise be dictated as another identity among the others that are listed.}, a direct computation reveals that in order for \eqref{hohoho} to satisfy the Jacobi identity 
\begin{equation} \label{gabon-B}
\begin{split}
	& [\rho(Y_3,X_2),X_1]_\mathcal{A} - [\rho(Y_3,X_1),X_2]_\mathcal{A} \\
	& \qquad \qquad = \rho(\sigma(X_2,Y_3),X_1)  - \rho(Y_3,[X_1,X_2]_\mathcal{A}) - \rho(\sigma(X_1,Y_3),X_2),  \\
	& \rho(Y_3,\rho(Y_2,X_1)) + \rho(\theta(Y_2,Y_3),X_1) - \rho(Y_2,\rho(Y_3,X_1))  \\
	& \qquad \qquad = -\psi(\sigma(X_1,Y_2),Y_3) - [\psi(Y_2,Y_3),X_1]_\mathcal{A} + \psi(\sigma(X_1,Y_3),Y_2) ,  \\
	& -\sigma(\rho(Y_2,X_1),Y_3) -\sigma(X_1,\theta(Y_2,Y_3)) + \sigma(\rho(Y_3,X_1),Y_2) \\
	& \qquad \qquad = -\theta(\sigma(X_1,Y_2),Y_3) + \theta(\sigma(X_1,Y_3),Y_2).
\end{split}
\end{equation}
and
\begin{equation} \label{gabon-A}
\begin{split}
		& \circlearrowright \psi(\theta(Y_1,Y_2),Y_3) - \circlearrowright\rho(Y_3,\psi(Y_1,Y_2)) = 0, \\
	& \circlearrowright \theta(\theta(Y_1,Y_2),Y_3) + \circlearrowright\sigma(\psi(Y_1,Y_2),Y_3) = 0 
 \end{split}
\end{equation}
are to be satisfied for any $X_1,X_2,X_3 \in \Gamma(\mathcal{A})$ and any $Y_1,Y_2,Y_3 \in \Gamma(\mathcal{B})$. Similarly, the Leibniz identity is satisfied if
\begin{equation} \label{gabon2}
\begin{split} 
	& \psi(Y_1,fY_2) = f\psi(Y_1,Y_2), \\
	& \theta(Y_1,fY_2) = f\theta(Y_1,Y_2) + \mathcal{L}_{\mathfrak{a}_{\mathcal{B}}(Y_1)}(f)Y_2 , \\
	& \rho(fY,X) = f\rho(Y,X),\\
	& \rho(Y,fX) = f\rho(Y,X) + \mathcal{L}_{\mathfrak{a}_{\mathcal{B}}(Y)}(f)X .
\end{split}
\end{equation}
Finally,
\begin{equation} \label{gabon3}
\begin{split}
	& [\mathfrak{a}_{\mathcal{B}}(Y),\mathfrak{a}_{\mathcal{A}}(X)] = \mathfrak{a}_{\mathcal{A}}(\rho(Y,X)) - \mathfrak{a}_{\mathcal{B}}(\sigma(X,Y)), \\
	& [\mathfrak{a}_{\mathcal{B}}(Y_1),\mathfrak{a}_{\mathcal{B}}(Y_2)] = \mathfrak{a}_{\mathcal{A}}(\psi(Y_1,Y_2)) + \mathfrak{a}_{\mathcal{B}}(\theta(Y_1,Y_2))
\end{split}
\end{equation}
are the compatibility conditions for the anchor map \eqref{bobobo}.

If the conditions in \eqref{gabon-B}-\eqref{gabon3} are satisfied, then the product bundle $\mathcal{A}\bowtie_\psi \mathcal{B}:=\mathcal{A}\times \mathcal{B}$ is a Lie algebroid, called the \emph{unified product} of $\mathcal{A}$ and $\mathcal{B}$. In accordance with the terminology of \cite{esen2021bicocycle}, we shall also call this product Lie algebroid the \emph{cocycle double cross product} of $\mathcal{A}$ and $\mathcal{B}$. These observations are summarized in the following theorem.

\begin{theorem}\label{thm-ula}
Let $(\mathcal{A},\tau,M,\mathfrak{a}_{\mathcal{A}},[\bullet,\bullet]_{\mathcal{A}})$ be a Lie algebroid, and $(\mathcal{B},\kappa,M)$ a vector bundle. Equipped with the maps \eqref{muratcan} and \eqref{rolar}, the product bundle 
\begin{equation}
(\mathcal{A}\bowtie_\psi \mathcal{B} := \mathcal{A}\times \mathcal{B},\tau_{\bowtie_\psi},M,\mathfrak{a}_{\bowtie_\psi},[\bullet,\bullet]_{\bowtie_\psi})
\end{equation} 
is a Lie algebroid along with the Lie bracket \eqref{hohoho}  and the anchor map  \eqref{bobobo} if and only if the conditions \eqref{gabon-B}-\eqref{gabon3} are satisfied. 
\end{theorem}

On the contrary, if a Lie algebroid is a product of a Lie subalgebroid of it and a complementary subbundle, then its structure is determined as a unified product. 

\begin{theorem}\label{mezgi}
    Let $(\mathcal{M},\upsilon,M, \mathfrak{a}, [\bullet,\bullet]_\mathcal{M} )$ be a Lie algebroid and 
    	\begin{equation}
    \Big(\mathcal{A},\tau=\upsilon\big \vert_\mathcal{A},M,\mathfrak{a}_\mathcal{A}=\mathfrak{a}\big \vert_\mathcal{A},[\bullet,\bullet] _\mathcal{A} = [\bullet,\bullet]_\mathcal{M}\big \vert_\mathcal{A} \Big)
	\end{equation}
 be a Lie subalgebroid so that $\mathcal{M} \cong \mathcal{A} \times \mathcal{B}$ as vector bundles, for some complementary vector bundle
    	\begin{equation}
    (\mathcal{B},\kappa=\upsilon\big \vert_\mathcal{B},M). 
     	\end{equation}
 Then $\mathcal{M}$ is a unified product Lie algebroid that is, $\mathcal{M}=\mathcal{A}\bowtie _ \psi \mathcal{B}$. In this case, the mappings in \eqref{muratcan} and \eqref{rolar} may be obtained from 
 	\begin{equation} \label{ezgi}
		[Y_1,Y_2] _{\mathcal{M}}   = \psi(Y_1,Y_2) + \theta(Y_1,Y_2),  \qquad 
  	[Y,X]_{\mathcal{M}} =  \rho(Y,X) - \sigma(X,Y),
     	\end{equation}
and they satisfy \eqref{gabon-A}-\eqref{gabon3}.
\end{theorem}

Let us conclude with the presentation in matrix components. To this end, let $(x^i)$ denote a local coordinate system on $M$, along with $(x^i,y^\alpha)$ on $\mathcal{A}$, and let $(e_\alpha)$ be the (projective) basis of $\Gamma(\mathcal{A})$. Let, similarly, $(e_a)$  represents a projective basis of $\Gamma(\mathcal{B})$, and let $(x^i,y^a)$ stand for a system of local coordinates on $\mathcal{B}$. Accordingly, by a slight abuse of notation, the projective basis for  $\Gamma(\mathcal{A}\bowtie _ \psi \mathcal{B}) \cong \Gamma(\mathcal{A})\times \Gamma(\mathcal{B})$ may be denoted by $(e_\alpha,e_a)$, while a set of local coordinates on which may be represented by $(x^i,y^\alpha,y^a)$.

Accordingly, the anchor map may be represented by 
	\begin{equation} \label{monon1} \mathfrak{a}_{\bowtie _ \psi }(e_\alpha) = (\mathfrak{a}_{\mathcal{A}})^i_\alpha\frac{\partial}{\partial x^i}, \qquad  
 \mathfrak{a}_{\bowtie _ \psi }(e_a) = (\mathfrak{a}_{\mathcal{B}})^i_a\frac{\partial}{\partial x^i},
	\end{equation}
while the representations in \eqref{rolar} may be expressed as 
	\begin{equation} \label{monon2}   \rho(e_a,e_\alpha)= \mathfrak{R}^{\beta}_{a\alpha } e_\beta,\qquad
 \sigma(e_\alpha,e_a)= - \mathfrak{S}^{b}_{a\alpha} e_b, 
	\end{equation}
 respectively. Finally, setting
	\begin{equation} \label{monon3} 
 \psi(e_a,e_b)=\mathfrak{P}^{\alpha}_{ab} e_\alpha,\qquad \theta(e_a,e_b)=\mathfrak{T}^{d}_{ab} e_d,
	\end{equation}
for the mappings in \eqref{muratcan}, the structure of the (unified) product Lie algebroid may be presented as
\begin{equation} \label{monon4}  
[e_\beta,e_\gamma]_{\bowtie _ \psi }=\mathfrak{H}_{\beta\gamma}^\alpha e_\alpha, \qquad  [e_b,e_d]_{\bowtie _ \psi }= \mathfrak{P}_{b d}^\alpha e_\alpha + \mathfrak{T}_{b d}^a e_a, \qquad [e_a,e_\alpha]_{\bowtie _ \psi } = \mathfrak{R}^{\beta}_{a\alpha } e_\beta + \mathfrak{S}_{a\alpha}^b e_b,
\end{equation}
where the former bracket uses the structure functions of the Lie algebroid $\mathcal{A}$, which were given above in \eqref{loc-algebroid}.

\subsubsection*{Unified Product Lie Algebras.} 

We find it instructive to record the particular case of the base manifold $M$ to be a single point. Along the lines of the above constructions, this particular choice corresponds to build a Lie algebra out of a Lie subalgebra and a complementary vector space, which has been studied in \cite{Agore-Lie}. We shall summarize the construction below as a corollary to Theorem \ref{thm-ula} and Theorem \ref{mezgi} above.

\begin{corollary}
Let $(\mathfrak{g},[\bullet,\bullet]_{\mathfrak{g}})$ be a Lie algebra, and let $\mathfrak{h}$ be a vector space, equipped with two sets of mappings
\begin{equation}\label{set-1}
\psi:\mathfrak{h} \times \mathfrak{h} \to \mathfrak{g}, \qquad \theta:\mathfrak{h} \times \mathfrak{h} \to \mathfrak{h}
\end{equation}
and
\begin{equation}\label{set-2}
\rho:\mathfrak{h} \times \mathfrak{g} \to \mathfrak{g}, \qquad \sigma:\mathfrak{g} \times \mathfrak{h} \to \mathfrak{h}. 
\end{equation}
Then $\mathfrak{g}\bowtie_\psi \mathfrak{h} := \mathfrak{g} \oplus \mathfrak{h}$ is a Lie algebra through
\begin{align*}
[(X_1,Y_1),(X_2,Y_2)]_{\bowtie_\psi} & := \Big([X_1,X_2]_ \mathcal{A}- \rho(Y_2,X_1) + \rho(Y_1,X_2) + \psi(Y_1,Y_2) ,\\
&\qquad \qquad \theta(Y_1,Y_2) - \sigma(X_2,Y_1) + \sigma(X_1,Y_2)\Big)
\end{align*}
for any $X_1,X_2 \in \mathfrak{g}$ and any $Y_1,Y_2 \in \mathfrak{h}$ if and only if (assuming the latter map in \eqref{set-2} is a Lie algebra action) the conditions listed in \eqref{gabon-B} and \eqref{gabon-A} are satisfied. Conversely, if $(\mathfrak{k}, [\bullet,\bullet]_\mathfrak{k})$ is a Lie algebra and $\mathfrak{g} \subseteq \mathfrak{k}$ is a Lie subalgebra so that $\mathfrak{k}\cong \mathfrak{g}\oplus\mathfrak{h}$ as vector spaces, for some complementary vector space $\mathfrak{h}$, then $\mathfrak{k}$ is a unified product of $\mathfrak{g}$ and $\mathfrak{h}$. The structure maps \eqref{set-1} and \eqref{set-2}, in this case, may be obtained by 
\[
[Y_1,Y_2] _{\mathfrak{k}}   = \psi(Y_1,Y_2) + \theta(Y_1,Y_2),  \qquad 
[Y,X]_{\mathfrak{k}} =  \rho(Y,X) - \sigma(X,Y).
\]
\end{corollary}

\subsection{Bicocycle Double Cross Product Lie Algebroids}\label{sec-bicocycle}

In the previous subsection, we presented the algebraic background of decomposing a Lie algebroid into a Lie subalgebroid and a vector bundle. In the present subsection, we shall upgrade that construction to the most general (and hence the most flexible) way of decomposing a Lie algebroid. Namely, we shall address algebraic foundations of decomposing a Lie algebroid into two vector bundles admitting all possible interactions. 

Let $(\mathcal{A},\tau,M)$ and $(\mathcal{B},\kappa,M)$ be two vector bundles along with the anchors $\mathfrak{a}_{\mathcal{A}}:\mathcal{A}\mapsto TM$ and $\mathfrak{a}_{\mathcal{B}}:\mathcal{B}\mapsto TM$, respectively. To the sets \eqref{muratcan} and \eqref{rolar} of skew-symmetric mappings we shall now adjoin another set
\begin{equation}\label{muratcan-ph}
\phi:\Gamma(\mathcal{A})\times \Gamma(\mathcal{A}) \longrightarrow \Gamma(\mathcal{A}), \qquad \zeta:\Gamma(\mathcal{A})\times \Gamma(\mathcal{A}) \longrightarrow \Gamma(\mathcal{B})
\end{equation}
of (also skew-symmetric) mappings.
In this case, a skew-symmetric bilinear operation on the space of sections of the product bundle $\mathcal{A}\times \mathcal{B}$, may be formulated by 	
\begin{equation}\label{hohoho+}
	\begin{split}
 	\Gamma(\mathcal{A}  \times  \mathcal{B}) \times  	\Gamma(\mathcal{A}  \times \mathcal{B}) & \longrightarrow 	\Gamma(\mathcal{A} \times  \mathcal{B}), \\		
		[(X_1,Y_1),(X_2,Y_2)]_{_\zeta\bowtie_\psi} & = \Big(\phi(X_1,X_2)- \rho(Y_2,X_1) + \rho(Y_1,X_2) + \psi(Y_1,Y_2) ,\\
		&\qquad \qquad \theta(Y_1,Y_2) - \sigma(X_2,Y_1) + \sigma(X_1,Y_2)+\zeta(X_1,X_2)\Big).
	\end{split}
\end{equation}
We wish to emphasize the mapping $\phi$ in the first entry of the bracket and the mapping  $\zeta$ in the second entry.
The product bundle may further be furnished with the bundle map
\begin{equation}\label{bobobo+}
\mathfrak{a}_{_\zeta\bowtie_\psi}:=\mathcal{A}\times \mathcal{B} \longrightarrow TM,\qquad 	 (X,Y) \mapsto \mathfrak{a}_{\mathcal{A}}(X) + \mathfrak{a}_{\mathcal{B}}(Y)
\end{equation} 
which shall serve as the anchor.

Tedious, though straightforward, calculations reveal at once that \eqref{hohoho+} satisfies the Jacobi identity if (and only if)
\begin{equation}\label{muratcanas1}
   \begin{split}
        &\rho(Y_3, \phi(X_1,X_2)) - \psi(\zeta(X_1,X_2), Y_3) + \phi(\rho(Y_3,X_2), X_1) \\
    & \qquad \qquad = \rho(\sigma(X_2,Y_3),X_1) + \phi(\rho(Y_3,X_1),X_2) - \rho(\sigma(X_1,Y_3),X_2) ,\\
    & \sigma(\phi(X_1,X_2),Y_3) + \theta(\zeta(X_1,X_2),Y_3) - \zeta(\rho(Y_3,X_2), X_1) \\
    & \qquad \qquad = \sigma(X_1,\sigma(X_2,Y_3)) - \zeta(\rho(Y_3,X_1),X_2) - \sigma(X_2,\sigma(X_1,Y_3)), \\
    & \rho(Y_3,\rho(Y_2,X_1)) + \psi(\sigma(X_1,Y_2),Y_3) + \phi(\psi(Y_2,Y_3),X_1) ,\\
    & \qquad \qquad = - \rho(\theta(Y_2,Y_3),X_1) + \rho(Y_2,\rho(Y_3,X_1)) + \psi(\sigma(X_1,Y_3),Y_2) ,\\
    & \sigma(\rho(Y_2,X_1),Y_3) - \theta(\sigma(X_1,Y_2),Y_3) - \zeta(\psi(Y_2,Y_3),X_1), \\&
     \qquad \qquad =  -\sigma(X_1,\theta(Y_2,Y_3)) + \sigma(\rho(Y_3,X_1),Y_2) - \theta(\sigma(X_1,Y_3),Y_2)\\
   \end{split}
\end{equation}
and
\begin{equation} \label{muratcanas2}
\begin{split}
		& \circlearrowright \phi(\phi(X_1,X_2),X_3) + \circlearrowright\rho(\zeta(X_1,X_2),X_3) = 0, \\
	& \circlearrowright \zeta(\phi(X_1,X_2),X_3) - \circlearrowright\sigma(X_3,\zeta(X_1,X_2)) = 0, \\
 & \circlearrowright \psi(\theta(Y_1,Y_2),Y_3) - \circlearrowright\rho(Y_3,\psi(Y_1,Y_2)) = 0, \\
	& \circlearrowright \theta(\theta(Y_1,Y_2),Y_3) + \circlearrowright\sigma(\psi(Y_1,Y_2),Y_3) = 0
 \end{split}
\end{equation}
are satisfied. 

Similarly, the Leibniz identity is satisfied if (and only if)
\begin{equation}\label{muratcanas4}
    \begin{split}
        & \phi(X_1,fX_2) = f\phi(X_1,X_2) + \mathcal{L}_{\mathfrak{a}_\mathcal{A}(X_1)}(f)X_2 ,\\
        &\zeta(X_1,fX_2) = f\zeta(X_1,X_2) ,\\
        & \rho(fY,X) = f\rho(Y,X) ,\\
        & \sigma(X,fY) = f \sigma(X,Y) + \mathcal{L}_{\mathfrak{a}_\mathcal{A}(X)}(f)Y ,\\
        & \rho(Y,fX) = f\rho(Y,X) + \mathcal{L}_{\mathfrak{a}_\mathcal{B}(Y)}(f)X ,\\
        & \sigma(fX,Y) = f\sigma(X,Y) ,\\
        & \psi(Y_1,fY_2) = f\psi(Y_1,Y_2), \\
        & \theta(Y_1,fY_2) = f\theta(Y_1,Y_2) + \mathcal{L}_{\mathfrak{a}_\mathcal{B}(Y_1)}(f)Y_2.
    \end{split}
\end{equation}
Finally \eqref{bobobo+} satisfies the anchor map property if
\begin{equation}\label{muratcanas3}
    \begin{split}
      & [\mathfrak{a}_{\mathcal{A}}(X_1),\mathfrak{a}_{\mathcal{A}}(X_2)] = \mathfrak{a}_{\mathcal{A}}(\phi(X_1,X_2)) + \mathfrak{a}_{\mathcal{B}}(\zeta(X_1,X_2)), \\
      & [\mathfrak{a}_{\mathcal{B}}(Y),\mathfrak{a}_{\mathcal{A}}(X)] = \mathfrak{a}_{\mathcal{A}}(\rho(Y,X)) - \mathfrak{a}_{\mathcal{B}}(\sigma(X,Y)) ,\\
       & [\mathfrak{a}_{\mathcal{B}}(Y_1),\mathfrak{a}_{\mathcal{B}}(Y_2)] = \mathfrak{a}_{\mathcal{A}}(\psi(Y_1,Y_2)) + \mathfrak{a}_{\mathcal{B}}(\theta(Y_1,Y_2)).
    \end{split}
\end{equation}

In accordance with the structure of the previous subsection, let us record these results in the form of a theorem.

\begin{theorem}\label{thm-muratcan} 
Given two vector bundles  $(\mathcal{A},\tau,M)$ and  $(\mathcal{B},\kappa,M)$, equipped with (skew-symmetric bilinear) mappings \eqref{muratcan},  \eqref{rolar}, and \eqref{muratcan-ph}, the product bundle $\mathcal{A} \times \mathcal{B}$ has the structure of a Lie algebroid through the Lie bracket \eqref{hohoho+}  and the anchor map  \eqref{bobobo+} if and only if the conditions listed in \eqref{muratcanas1}-\eqref{muratcanas3} are all satisfied. 
\end{theorem}

Following the terminology of \cite{esen2021bicocycle}, in order to highlight the \emph{twisted cocycles} we shall denote the product Lie algebroid as $\mathcal{A} {\hspace{.1cm}}_{\zeta\hspace{-.1cm}} \bowtie_\psi \mathcal{B} := \mathcal{A} \times \mathcal{B}$, and call it the \emph{bicocycle double cross product} of $\mathcal{A}$ and $\mathcal{B}$.

This theorem generalizes the unified product Lie algebroid structure to the sum of two vector bundles. More concretely, as a particular instance, if $\zeta$ in \eqref{muratcan-ph} is identically zero, then one arrives at Theorem \ref{thm-ula} where unified product Lie algebroid structure is exhibited.  

In contrast to the complexity of the conditions \eqref{muratcanas1}-\eqref{muratcanas3}, the bicocycle double cross product construction offers utmost flexibility when it comes to decompose a Lie algebroid. More precisely, we have the following.

\begin{theorem}\label{mezgi+}
Given a Lie algebroid $(\mathcal{M},\upsilon,M,[\bullet,\bullet]_\mathcal{M}, \mathfrak{a} )$, let $\mathcal{M} \cong \mathcal{A} \times \mathcal{B}$, as vector bundles, for two complementary vector bundles
    	\begin{equation}
     (\mathcal{A},\tau=\upsilon\big \vert_\mathcal{A},M  ),\qquad  (\mathcal{B},\kappa=\upsilon\big \vert_\mathcal{B},M). 
	\end{equation}
Then $\mathcal{M}$ is the bicocycle double cross product of $\mathcal{A}$ and $\mathcal{B}$, and the structure maps are obtained by   
\begin{equation} \label{ezg}
  \begin{split}
  	[X_1,X_2] _{\mathcal{M}}   & = \phi(X_1,X_2) + \zeta(X_1,X_2),  \\
      		[Y_1,Y_2] _{\mathcal{M}}  & = \psi(Y_1,Y_2) + \theta(Y_1,Y_2),  \\
  	[Y,X]_{\mathcal{M}} &=  \rho(Y,X) - \sigma(X,Y) .
  \end{split}
\end{equation}
\end{theorem}

Let us conclude, as usual, with the coordinates in the finite dimensional setting. To this end, let $(x^i)$ be a local coordinate system on $M$, and let $(e_\alpha)$ and $(e_a)$ be projective bases of $\Gamma(\mathcal{A})$ and $\Gamma(\mathcal{B})$ respectively. 
In addition to \eqref{monon2} and \eqref{monon3}, we now have
\begin{equation} \label{monon5} 
 \phi(e_\alpha,e_\beta)=\mathfrak{H}^{\gamma}_{\alpha\beta} e_\gamma,\qquad \zeta(e_\alpha,e_\beta)=\mathfrak{Z}^{a}_{\alpha\beta} e_a,
\end{equation}
where $\mathfrak{H}_{\beta\gamma}^\alpha$ now fail to be proper structure functions due to the existence of the twisted cocycle terms $\mathfrak{Z}^{a}_{\alpha\beta}$. Accordingly the Lie bracket is determined by
\begin{equation} \label{monon6}  
	  [e_\beta,e_\gamma]_{_\zeta\bowtie _ \psi }=\mathfrak{H}_{\beta\gamma}^\alpha e_\alpha + \mathfrak{Z}^{a}_{\beta\gamma} e_a, \qquad  [e_b,e_d]_{_\zeta\bowtie _ \psi }= \mathfrak{P}_{b d}^\alpha e_\alpha + \mathfrak{T}_{b d}^a e_a, \qquad [e_a,e_\alpha]_{_\zeta\bowtie _ \psi } = \mathfrak{R}^{\beta}_{a\alpha } e_\beta + \mathfrak{S}_{a\alpha}^b e_b,
	\end{equation}
while the anchor map is computed to be 
\begin{equation} \label{monon11} \mathfrak{a}_{_\zeta\bowtie _ \psi }(e_\alpha) = (\mathfrak{a}_{\mathcal{A}})^i_\alpha\frac{\partial}{\partial x^i}, \qquad  
 \mathfrak{a}_{_\zeta\bowtie _ \psi }(e_a) = (\mathfrak{a}_{\mathcal{B}})^i_a\frac{\partial}{\partial x^i}.
	\end{equation}

\subsubsection*{Bicocycle Double Cross Product Lie Algebras.}

In the particular case of $\mathcal{M}$ being a Lie algebroid over a point, the space $\Gamma(\mathcal{M})$ of sections of $\mathcal{M}$ happens to be a Lie algebra. Then, \cite[Prop. 2.1]{esen2021bicocycle} and \cite[Prop. 2.2]{esen2021bicocycle} follow from Theorem \ref{thm-muratcan} and Theorem \ref{mezgi+} respectively.

\begin{corollary}
Let $(\mathfrak{g},\mathfrak{h})$ be a pair of vector spaces equipped with mappings \eqref{muratcan},  \eqref{rolar} and \eqref{muratcan-ph} (considering $\mathfrak{g} = \Gamma(\mathcal{A})$ and $\mathfrak{h} = \Gamma(\mathcal{B})$). Then, the direct sum vector space $\mathfrak{g} {\hspace{.1cm}}_{\zeta\hspace{-.1cm}} \bowtie_\psi \mathfrak{h} := \mathfrak{g} \oplus \mathfrak{h}$ is a Lie algebra by the bracket \eqref{hohoho+} if and only if the conditions \eqref{muratcanas2} and \eqref{muratcanas1} are satisfied. Conversely, if a Lie algebra  $(\mathfrak{k},[\bullet,\bullet]_\mathfrak{k})$ can be realized as a direct sum $\mathfrak{k} \cong \mathfrak{g} \oplus \mathfrak{h}$ of complementary vector spaces $\mathfrak{g}$ and $\mathfrak{h}$, then $\mathfrak{k}=\mathfrak{g}{\hspace{.1cm}}_{\zeta\hspace{-.1cm}} \bowtie_\psi \mathfrak{h}$ as Lie algebras.  
\end{corollary}

\subsection{Twisted Poisson Algebroids and Dirac Structures} \label{Dirac-sec}

A manifold $P$ is called an \emph{almost Poisson manifold} if the  space $\mathcal{C}^\infty(P)$ of smooth functions is equipped with a bilinear, skew-symmetric bracket $\{\bullet,\bullet\}$ which satisfy the Leibniz identity
\begin{equation}
\{F_1,F_2 F_3\} = \{F_1,F_2\}F_3 + F_2\{F_1,F_3\}, 
\end{equation}
for any $F_1,F_2,F_3 \in \mathcal{C}^\infty(P)$. Accordingly, the bracket $\{\bullet,\bullet\}$ may be represented by a bivector field $\Lambda$ as 
\begin{equation} \label{bivec-PoissonBra}
\Lambda(dF_1,dF_2) = \{F_1,F_2\},
\end{equation}
where $dF_1,dF_2$ denote the de-Rham exterior derivatives. An almost Poisson manifold can thus be represented by a pair $(P,\Lambda)$ of a manifold and a bivector field.

For a Hamiltonian function $\mathcal{H}$ defined on an almost Poisson manifold $P$, the Hamilton's equation is defined as
\begin{equation}
\frac{du}{dt} = \{u,\mathcal{H}\},
\end{equation}
where $u$ is in $P$. Accordingly, the Hamiltonian vector field $X_{\mathcal{H}}$ for a Hamilton function $\mathcal{H}$ is defined to be
\begin{equation}
X_{\mathcal{H}}(F) = \{F,\mathcal{H}\}.
\end{equation}
If there exists a non-constant function $C$ such that for all $F$:
\begin{equation}
\{F,C\} = 0
\end{equation}
then the almost Poisson framework is called degenerate. Such a  function $C$ is called a  \emph{Casimir function}. It should be noted that a Casimir function cannot cause any dynamics since the Hamiltonian vector field is identically zero for a Casimir function. 

On an almost Poisson manifold $(P,\{\bullet,\bullet\})$ the mapping
\begin{equation} 
\begin{split}
    \mathfrak{J}: \wedge^3  \mathcal{C}^\infty(P) \longrightarrow \mathcal{C}^\infty(P),   \qquad
(F_1, F_2, F_3)\mapsto~ \circlearrowright \{F_1,  \{F_2, F_3\} \} .
\end{split}
\end{equation}
is called the \emph{Jacobiator}. If the Jacobiator vanishes, then the bracket on the almost Poisson manifold $P$ satisfies the Jacobi identity. In this case the almost Poisson manifold $(P,\{\bullet,\bullet\})$ is called a \emph{Poisson manifold}. 

From the point of view of bivector fields, the almost Poisson manifold $(P,\Lambda)$ to be a Poisson manifold corresponds to 
\begin{equation} \label{Poisson-cond}
[\Lambda,\Lambda] = 0,
\end{equation} 
the bracket on the left hand side being the Schouten-Nijenhuis bracket.

In this case, the characteristic distribution of the Poisson manifold $P$, that is the image space of all Hamiltonian vector fields, is integrable. Namely,
\begin{equation}
\label{eq-brackets+mc}
X_{\{F_1,F_2\}}+[X_{F_1},X_{F_2}]=0.\end{equation}
This implies also that the manifold $P$ is foliated by symplectic leaves. This is a consequence of the non-degeneracy of the Poisson framework on each leaf. 

Consider an almost Poisson manifold $(P,\Lambda)$. The musical mapping $\Lambda^\sharp$  induced by the bivector field $\Lambda$ is defined to be
\begin{equation} \label{sharp}
\Lambda^\sharp: \Gamma^1(P)\longrightarrow \mathfrak{X}(P), \qquad \langle \beta, \Lambda^\sharp(\alpha)\rangle = \Lambda(\alpha, \beta),
\end{equation}
where the pairing on the left-hand side is the one between the space $\Gamma^1(P)$ of one-form sections and the space $\mathfrak{X}(P)$ of vector fields. 
Then sections of the cotangent bundle $\pi_P:T^*P\mapsto P$ admit a skew-symmetric bilinear operation
\begin{equation} \label{alg-bra}
[\alpha ,\beta] =\mathfrak{L}_{\Lambda^\sharp(\alpha)}\left( \beta \right) -\mathfrak{L}_{\Lambda^\sharp(\beta) }\left( \alpha \right)
- d (\Lambda(\alpha,\beta) ).
\end{equation}
A direct observation reads that this bracket enables us to determine an \emph{almost Lie algebroid} $(T^*P,\pi_P, P,\Lambda^\sharp,[\bullet,\bullet])$ where the musical mapping $\Lambda^\sharp$ is the anchor map. We denote this algebroid by the following commutative diagram:
\begin{equation} 
\xymatrix{
T^*P \ar[rr]^{\Lambda^\sharp} \ar[dr]_{\pi_P} &  & TP \ar[dl]^{\tau_P} \\
& P
}
\end{equation}
If $\Lambda$ is a Poisson bivector field that satisfies the Jacobi identity \eqref{Poisson-cond} then the Jacobiator for the bracket \eqref{alg-bra} vanishes identically either. Accordingly, the bracket \eqref{alg-bra} satisfies the Jacobi identity, and the quintuple 
\begin{equation} 
(T^*P,\pi_P, P,\Lambda^\sharp,[\bullet,\bullet])
\end{equation}
turns out to be a genuine Lie algebroid.

\subsubsection*{Dirac Bundle Realization.}

An \emph{almost Dirac structure}\footnote{It is common, in the literature, to represent the Dirac structure as a subbundle of $TM\times T^*M$. In the present manuscript, however, we shall reverse the order of the bundles, and determine the Dirac structure as a subbundle of Whitney sum $T^*M\times TM$.} on a manifold $M$ is a subbundle $D$ of the Whitney sum $T^*M\times TM$ such that $D$ is a maximal isotropic subbundle of $T^*M\times TM$ with respect to the pairing \cite{bursztyn2013brief, courant1990dirac, yoshimura2006dirac3,yoshimura2006dirac1, yoshimura2006dirac2}
\begin{equation}\label{pairing}
(T^*M\times TM) \oplus (T^*M\times TM) \longrightarrow \mathbb{R},\qquad 
(v_1,u_1) \cdot (v_2,u_2) =u_1(v_2)+u_2(v_1).
\end{equation} 
Therefore we denote a Dirac structure by
\begin{equation} 
D\subset T^*M\times TM,\qquad  D=D^\perp.
\end{equation} 
Define a bundle structure $\tau:D \mapsto M$ by projecting a pair $(v,u)$ to the base point in $M$. 
If $pr_2(D)$ has a constant rank, then we say that $D$ is a regular almost Dirac structure, and $pr_2(D)$ defines a regular foliation.  

An almost Dirac structure is called a Dirac structure if the space $\Gamma(D)$  of sections for the fibration $D\mapsto M$ is closed with respect to the (skew-symmetric)  \emph{Courant bracket}
\begin{equation}\label{courantbracket}
\begin{split}
   [\bullet,\bullet]&:\Gamma(D)\oplus\Gamma(D)\longrightarrow \Gamma(D), \\ 
   &\qquad 
[(\alpha_1, U_1),(\alpha_2, U_2)]] = (\mathcal{L}_{U_1}\alpha_2 -\mathcal{L}_{U_2}\alpha_1 +\frac{1}{2}d (\alpha_1(U_2)-\alpha_2(U_1)), [U_1, U_2 ]). 
\end{split}
\end{equation}
A Dirac structure
$D$ on a manifold $M$ carries a Lie algebroid structure with the anchor map  $pr_2: D\mapsto TM$ (projection to the second factor in $D$). The Lie algebroid $D$ is determined by
\begin{equation} 
\xymatrix{
D \ar[rr]^{pr_2} \ar[dr]_{\tau} &  & TM \ar[dl]^{\tau_M} \\
& M
}
\end{equation}
and we shall denote it by
\begin{equation} 
D=(D,\tau, M,pr_2,[\bullet,\bullet]).
\end{equation}

Let us now examine the Lie algebroid $D$ from the point of view of the unified products. We remark that the Dirac bundle $D$ decomposes into the direct sum of two of its subbundles. The decomposition may simply be given by 
\begin{equation} \label{Diracdecomp}
(pr_1,pr_2)\circ \Delta :D\longrightarrow D^{T^*}\times D^{T},
\end{equation}
where $\Delta$ is the diagonal map that takes an element in $D$ to a two copies of it in $D\times D$, while $D^T$ and $D^{T^*}$ stands for the tangent and the cotangent bundle parts of $D$ respectively. It is immediate to observe  that there is no cocycle term ($\psi=0$) reading that both $D^T$ and $D^{T^*}$ are Lie subalgebroids of $D$. So that we are in the realm of double cross product Lie algebroid category. In particular, the mutual actions are computed through 
\begin{equation} 
[U,\alpha] = \Big(\mathcal{L}_{U}\alpha - \frac{1}{2}d \alpha(U), 0 \Big)
\end{equation}
which gives that 
\begin{equation} \label{Dirac-mutu}
\rho(U,\alpha)= \mathcal{L}_{U}\alpha - \frac{1}{2}d \alpha(U),\qquad \sigma(\alpha,U)=0 .
\end{equation}
Evidently this reads a semi-direct product Lie algebroid decomposition
\begin{equation} 
D=D^{T^*}\rtimes D^T
\end{equation}
as a particular instance of double cross product Lie algebroid structure where the action term $\sigma$ is identically zero. 

\subsubsection*{The Twisted Category.} 

There is an interesting subcategory of almost Poisson manifolds.  In this case, the Jacobi identity \eqref{Poisson-cond} is still not satisfied, but it is equal to the pull-back of a closed three-form \cite{KlSt02,SeWe01}. Let us depict this geometry in more detail. 
It is possible to generalize the musical mapping $\Lambda^\sharp$ in \eqref{sharp} to the space $\Gamma^k(P)$ of $k$-forms on $P$ as well. We use the same notation for this mapping hoping no confusion may arise. We define
\begin{equation}\label{sharp-gen}
\Lambda^\sharp:\Gamma^k(P)\longrightarrow \mathfrak{X}^k(P), \qquad \langle \alpha_1\wedge\dots \wedge \alpha_k, \Lambda^\sharp(\omega) \rangle =(-1)^k\langle \omega, \Lambda^\sharp(\alpha_1)\wedge\dots \wedge \Lambda^\sharp(\alpha_k)\rangle, 
\end{equation} 
where the notation $\Lambda^\sharp$ on the right-hand side of the equality is the musical mapping \eqref{sharp} between the first-order tensor fields. By loosening the Jacobi identity \eqref{Poisson-cond} as
\begin{equation} \label{twisted-Poisson-cond}
\frac{1}{2}[\Lambda,\Lambda]=\Lambda^\sharp(\varphi)
\end{equation} 
for a closed three-form $\varphi$, one defines twisted Poisson structure. Here, $\Lambda^\sharp$ notation refers to the mapping in \eqref{sharp-gen} for $k=3$. We denote a twisted Poisson manifold by a three-tuple $(P,\Lambda,\varphi)$. In terms of the functions, the condition \eqref{twisted-Poisson-cond} is 
\begin{equation}  \label{twisted-Poisson-cond-bra}
\circlearrowright\{F_1,\{F_2,F_3\}\}=\varphi(\Lambda^\sharp(dF_1),\Lambda^\sharp(dF_2),\Lambda^\sharp(dF_3)),
\end{equation}
where $\Lambda^\sharp$ is the musical mapping in \eqref{sharp}.  
Finally, we note that
\begin{equation}
\label{eq-brackets}
X_{\{F_1,F_2\}}+[X_{F_1},X_{F_2}]=-\Lambda^\sharp(\varphi(X_{F_1},X_{F_2},\bullet)).\end{equation}

Recall the bracket \eqref{alg-bra} defined on the space of one-form sections over an almost Poisson manifold. Evidently, this bracket does not satisfy the Jacobi identity for the case of twisted Poisson manifolds. On the other hand, we can introduce another bracket 
\begin{equation} 
[\alpha, \beta]_\varphi=[\alpha, \beta]+\iota_{\Lambda^{\sharp}(\alpha) \wedge \Lambda^{\sharp}(\beta)}(\varphi),
\end{equation}
which satisfies the Jacobi identity. Further, the quintuple $\left(T^*P, \pi_P,P,\Lambda ^{\sharp},[\cdot, \cdot]_\varphi\right)$ is a Lie algebroid. 

Consider an almost Dirac structure $D$ on a manifold $M$ and, $\varphi$-twisted skew-symmetric Courant bracket \cite{courant1988beyond,SeWe01}
\begin{equation}
\label{eq-courant} 
\begin{split}
  [(\alpha_1,U_1),(\alpha_2,U_2)]_\varphi  &= [(\alpha_1,U_1),(\alpha_2,U_2)]+ \varphi(U_1,U_2,\bullet) \\
&=  (\mathcal{L}_{U_1}\alpha_2 -\mathcal{L}_{U_2}\alpha_1 + \frac{1}{2}d (\alpha_1(U_2)-\alpha_2(U_1))+ \varphi(U_1,U_2,\bullet), [U_1, U_2 ]) .
\end{split}
\end{equation}
An almost Dirac structure $D$ of $T^*M\times TM$ is a $\varphi$-twisted Dirac structure 
if it satisfies the integrability condition with respect to the $\varphi$-twisted skew-symmetric Courant bracket
$[\bullet, \bullet]_{\varphi}$.
As in the ordinary case, a twisted Dirac structure $D$ on $M$ induces a Lie algebroid on $D$ given by
the anchor map $pr_2|_{D}$ and the bracket $[\bullet, \bullet]_{\varphi}$.  Therefore $pr_2(D)$ is an integrable distribution on $M$. If $D$ is a regular almost Dirac structure such that $pr_2(D)\subset TM$ is an integrable distribution on $M$, then, there exists an exact  
three-form $\varphi$ with respect to which $D$ is a $\varphi$-twisted Dirac structure \cite{NaranjoBalseiro}. We denote this Lie algebroid by
\begin{equation}
\label{eq-courant+} 
D_{\varphi}=(D,\tau, M,pr_2,[\bullet,\bullet]_{\varphi}).
\end{equation}

In this uniform product analysis, a more complicated decomposition of the Dirac bundle is encountered. Once again, referring to \eqref{Diracdecomp}, we write $D=D^{T^*}\times D^{T}$. Notice that $D^{T^*}$ is remaining to be a Lie subalgebroid while $D^T$ is not. Let us examine the interaction between these subbundles step by step. The mutual \textit{actions} of the bundles are exactly the same with those given in \eqref{Dirac-mutu}. If we take the bracket
\begin{equation} 
  [(0, U_1),(0, U_2)]_\varphi  = \Big(\varphi(U_1,U_2,\bullet), [U_1, U_2 ]\Big)
\end{equation}
 of two tangent elements, the mappings in \eqref{muratcan} are interpreted as
\begin{equation} 
  \theta(U_1,U_2) =[U_1, U_2 ] \in \Gamma(D^T),\qquad  \psi(U_1, U_2 )  = \varphi(U_1,U_2,\bullet) \in \Gamma(D^{T^*}), 
\end{equation}
that is, $\varphi(U_1,U_2,\bullet)$ serves as the twisted cocycle term. So we decompose this Dirac bundle as 
\begin{equation} 
D=D^{T^*}\bowtie_\psi D^{T} 
\end{equation}
since one of the actions is trivial.

\subsection{Decomposition of Jacobi Lie Algebroid} \label{decomsec}

 A manifold $M$ equipped with a vector field $\mathcal{E}$ and a bivector field $\Lambda$ is a \emph{Jacobi manifold} if 
\begin{equation}\label{ident-Jac}
        [\Lambda,\Lambda] = 2 \mathcal{E} \wedge \Lambda, \qquad 
       [\mathcal{E},\Lambda] = 0,
\end{equation}
    where $[\bullet,\bullet ]$ is the Schouten-Nijenhuis bracket, see, for example, \cite{Lichnerowicz-Poi,Lichnerowicz-Jacobi,Marle-Jacobi,vaisman2002jacobi}. We denote a Jacobi manifold by a triple $(M,\Lambda,\mathcal{E})$. Starting with a Jacobi manifold, one may define an antisymmetric bilinear bracket
   \begin{equation}\label{bra-Jac}
        \{F_1,F_2\} = \Lambda(dF_1, dF_2) + F_1\mathcal{E}(F_2) - F_2\mathcal{E} (F_1)
 \end{equation}
satisfying the Jacobi identity. Furthermore, it fulfills the so-called weak Leibniz identity
\begin{equation}
        \operatorname{supp}(\{F_1,F_2\}) \subseteq \operatorname{supp} (F_1) \cap \text{supp} (F_2).
\end{equation}
This observation reads that the algebra is a local Lie algebra in the sense of Kirillov \cite{Kirillov}. The inverse of this assertion is also true, that is, a local Lie algebra determines a Jacobi structure. Important examples of Jacobi structures include locally conformally symplectic structures \cite{Banyaga2002,Bazzoni2018,esen21lcs,vaisman}, locally conformally cosymplectic structures \cite{Atesli23,Chinea91},  contact manifolds \cite{Arn}, and cocontact manifolds \cite{atecsli2023non,Leon-coco,Gaset2323,Rivas23}.

Consider a Jacobi manifold determined by the triplet $(M, \Lambda, \mathcal{E})$. 
We determine a homomorphism from the space of one-forms to the space of vector fields. It is determined by
\begin{equation}
\Lambda^\sharp:\Gamma^{1}(M) \longrightarrow \mathfrak{X}(M),\qquad \langle \Lambda^\sharp(\mu), \nu \rangle =\Lambda(\mu, \nu)
\end{equation}
for all $\mu$ and $\nu$ in $\Gamma^{1}(M)$. For a smooth  real-valued Hamiltonian function $\mathcal{H}$, the  Hamiltonian vector field $X_{\mathcal{H}}$ is defined by
\begin{equation}
X_{\mathcal{H}}=\Lambda^\sharp(d \mathcal{H})+\mathcal{H} \mathcal{E}.
\end{equation}
Note that, for the constant function $H=1$ the Hamiltonian vector field is $\mathcal{E}$. The mapping taking a Hamiltonian function $\mathcal{H}$ to the Hamiltonian vector field $X_\mathcal{H}$ is a Lie algebra homomorphism satisfying
\begin{equation}
\left[X_{F_1}, X_{F_2}\right]=X_{\{F_1, F_2\}}.
\end{equation}

The cotangent bundle of a Poisson manifold admits a Lie algebroid structure \cite{BhasVisw88,CostDazoWein87}. For a Jacobi manifold $(M, \Lambda, \mathcal{E})$, the picture is as follows. Consider the extended cotangent bundle $\mathbb{R} \times T^{*} M$ as the total space of the first jet prolongation of the fibration $M\mapsto \mathbb{R}$  \cite{Saunders-book}. The space of sections of $p_M:\mathbb{R} \times T^{*} M\mapsto M$ is the product space of real valued functions and one-forms $\mathcal{F}(M) \times \Gamma^1(M)$. So it consists of pairs $(F, \mu)$ where $F$ be a real valued function and $\mu$ is a one-form. A bracket on $\mathcal{F}(M) \times \Gamma^1(M)$ is defined to be
\begin{equation}\label{Jac-Brack}
\begin{split} 
\big\{(F_1, \mu),(F_2, \nu)\big\}_J&=\Big(\Lambda(\nu, \mu)+\Lambda^\sharp(\mu)(F_2)-\Lambda^\sharp(\nu)(F_1)+F_1 \mathcal{E}(F_2)-F_2 \mathcal{E}(F_1), \\
&\qquad \mathcal{L}_{\Lambda^\sharp(\mu)} \nu-\mathcal{L}_{\Lambda^\sharp(\nu)} \mu-d(\Lambda(\mu, \nu))+F_1 \mathcal{L}_{\mathcal{E}} \nu-F_2 \mathcal{L}_{\mathcal{E}} \mu-\iota_{\mathcal{E}}(\mu \wedge \nu) \Big),
\end{split}
\end{equation}
where $\Lambda^\sharp$ is the musical mapping induced by the bivector field $\Lambda$. Equipped with this bracket, the quintuple  $(\mathbb{R} \times T^{*} M,p_M,M,\mathcal{E}+\Lambda^\sharp,\{\bullet,\bullet\}_J)$ happens to be a Lie algebroid, \cite{KerbSoui93}. Here, the anchor map is given by 
\begin{equation}
\mathcal{E}+\Lambda^\sharp:\mathcal{F}(M) \times \Gamma^1(M) \longrightarrow \mathfrak{X}(M),\qquad (F,\mu)\mapsto \Lambda^\sharp(\mu)+F\mathcal{E}. 
\end{equation}
We refer the reader to \cite{LeonLopeMarrPadr03,LeonMarrJuan97,LeonMarr97} for further analysis on this Lie algebroid. If $\mathcal{E}$ is identically zero, then one can consider the cotangent bundle $T^*M$ as the total space and the bracket \eqref{Jac-Brack} reduces to the one on the algebra of one-forms on the Poisson manifold (hence the Poisson Lie algebroid structure). 

In view of Theorem \ref{mezgi+} we shall now examine this Lie algebroid from the point of view of the unified products. More precisely, we show that the Jacobi Lie algebroid fits into the unified product framework. 

To this end, we first determine the bracket on the function space by computing the bracket of two functions as 
\begin{equation}\label{Jac-Brack+-+-} 
\{F_1, F_2\}_J = \big(F_1 \mathcal{E}(F_2)-F_2 \mathcal{E}(F_1) , 0 \big).
\end{equation}
This reads the space $\mathcal{F}(M)$ of functions as a Lie subalgebroid of the Jacobi algebroid. Furthermore, referring to the previous calculation, we see the Lie bracket on the induced Lie subalgebroid as 
\begin{equation}
    \{F_1, F_2\}_{\mathcal{F}(M)}= F_1 \mathcal{E}(F_2)-F_2 \mathcal{E}(F_1).
\end{equation}
On the other hand, for the space of one-form sections, the situation is different. This space is not closed under the induced Lie bracket:  
\begin{equation}\label{Jac-Brack+-} 
\{\mu, \nu\}_J = \Big( \Lambda(\nu, \mu), \mathcal{L}_{\Lambda^\sharp(\mu)} \nu-\mathcal{L}_{\Lambda^\sharp(\nu)} \mu-d(\Lambda(\mu, \nu)) -\iota_{\mathcal{E}}(\mu \wedge \nu) \Big), 
\end{equation}
since there appears a term in the first slot manifesting the existence of twisted cocycle term. Accordingly, taking \eqref{ezgi} into account, we write  
\begin{equation}\label{Jac-Brack++} 
\begin{split} 
\theta (\mu,\nu) &=  \mathcal{L}_{\Lambda^\sharp(\mu)} \nu-\mathcal{L}_{\Lambda^\sharp(\nu)} \mu-d(\Lambda(\mu, \nu)) -\iota_{\mathcal{E}}(\mu \wedge \nu),  \\ 
   \psi (\mu,\nu) &=  \Lambda(\nu, \mu).
  \end{split}
\end{equation}
Notice that we have that the bivector field turns out to be the twisted cocycle term. 
Further, we compute the mutual actions using the bracket of cross terms as
\begin{equation}\label{Jac-Brack--}
\begin{split} 
\{\mu, F \}_J&=\Big(\Lambda^\sharp(\mu)(F),  -F \mathcal{L}_{\mathcal{E}} \mu \Big),
\end{split}
\end{equation}
which gives the following mappings: 
\begin{equation}
\rho (\mu,F)= \Lambda^\sharp(\mu)(F) ,\qquad \sigma(F,\mu)= F \mathcal{L}_{\mathcal{E}} \mu . 
\end{equation}
Therefore we conclude that the Jacobi Lie algebroid is a unified product of the Poisson Lie algebroid with the algebra of functions.

\section{Dynamics on Product Lie Algebroids} 
\label{dyndyn}

\subsection{Reversible Dynamics on Lie Algebroids}
\label{dyndyn-1}

Let $(\mathcal{A},\tau,M,\mathfrak{a}_{\mathcal{A}},[\bullet,\bullet]_{\mathcal{A}})$ denote a Lie algebroid. We express its dual bundle as $(\mathcal{A}^{\ast}, \pi, M)$, where $\pi$ stands for the projection dual to $\tau$. The fibers of $\mathcal{A}^{\ast}$ over each point $m$ in the base $M$, denoted as $\mathcal{A}^{\ast}_{m} = \pi^{-1}(m)$, are the linear algebraic dual vector spaces of the fibers of $\mathcal{A}$ corresponding to $\mathcal{A}_{m} = \tau^{-1}(m)$, see \cite{MackenzieDG,mackenzie1995lie,weinstein1996groupoids}.

The space of functions on the dual space $\mathcal{A}^{\ast}$ is spanned by two classes of functions. The first class consists of constant functions on fibers, i.e., functions on the base $M$. The second class consists of the linear functions on fibers, i.e., the sections of the projection $\tau$. Let us analyze the first class. Given a function $f$ on $M$, we can define a function $\widehat{f}$ on $\mathcal{A}^{\ast}$ as follows:
\begin{equation}\label{lolo}
\widehat{f}: \mathcal{A}^{\ast} \longrightarrow \mathbb{R}, \qquad \widehat{f}(u) = f \circ \pi(u).
\end{equation}
The second class will be obtained from the sections of the Lie algebroid $\mathcal{A}$. To this end, let us take a section $X$. The linear function $\widehat{X}$ on $\mathcal{A}^{\ast}$, for every $u$ in $\mathcal{A}^{\ast}$, is defined as:
\begin{equation}\label{nono}
\widehat{X}: \mathcal{A}^{\ast} \longrightarrow \mathbb{R}, \qquad \widehat{X}(u) = \langle u, X \circ \pi(u) \rangle.
\end{equation}
For these functions, a Poisson bracket on $\mathcal{A}^{\ast}$ is defined as:
\begin{equation} 
\{ \widehat{f}, \widehat{g} \}_{\mathcal{A}^{\ast}}  = 0, \qquad 
\{ \widehat{f}, \widehat{X} \}_{\mathcal{A}^{\ast}}  = \widehat{\mathcal{L}_{\mathfrak{a}_{\mathcal{A}} (X)} f}, \qquad 
\{ \widehat{X}_1, \widehat{X}_2 \}_{\mathcal{A}^{\ast}}  = \widehat{[X_1, X_2]}_{\mathcal{A} }. 
\label{Poisson on dual A}
\end{equation}
Here, the value of the second bracket on the right-hand side is the directional derivative of the function $f$ in the direction of the vector field $\mathfrak{a}_{\mathcal{A}}(X)$.  

 In finite dimensions, let $(x^i)$ be a local coordinate system on $M$ and $(x^i, y^\alpha)$ be a local coordinate system on $\mathcal{A}$. In parallel, choose the basis set $(e_\alpha)$ for the section space $\Gamma(\mathcal{A})$. 
Let us choose the dual basis $(e^\alpha)$ for the section space of $\mathcal{A}^{\ast}$. 
Assuming  $ (x^{i} )$, $ (x^{i},y^{\alpha } )$, and $ (x^{i},y_{\alpha } )$   local coordinates on $M$, $\mathcal{A}$, and $\mathcal{A}^{\ast}$ respectively, the Poisson brackets in \eqref{Poisson on dual A} take the form \cite{martinez2001,WeinLag,weinstein1996groupoids}:
\begin{equation} 
\{ x^{i}, x^{j} \}_{\mathcal{A}^{\ast}}  = 0, \qquad 
\{x^{i}, y_{\alpha} \}_{\mathcal{A}^{\ast}}  = (\mathfrak{a}_{\mathcal{A}})_{\alpha}^{i}, \qquad 
\{ y_{\alpha}, y_{\beta} \}_{\mathcal{A}^{\ast}}  = \mathfrak{H}_{\alpha \beta}^{\gamma} y_{\gamma}, 
\end{equation}
where $ (\mathfrak{a}_{\mathcal{A}})_{\alpha}^{i}$ represents the components of the anchor map in the chosen bases, and $ \mathfrak{H}_{\alpha \beta}^{\gamma} $ are the structure functions determined in equation \eqref{loc-algebroid}. In light of these local representations, the Poisson bivector corresponding to the Poisson bracket is computed to be
\begin{equation}\label{Poisson-A*}
\Lambda_{\mathcal{A}^{\ast}} = \frac{1}{2} \mathfrak{H}_{\alpha \beta}^{\gamma} y_{\gamma} \frac{\partial}{\partial y_{\alpha}} \wedge \frac{\partial}{\partial y_{\beta}} + (\mathfrak{a}_{\mathcal{A}})_{\alpha}^{i} \frac{\partial}{\partial x^{i}}\wedge \frac{\partial}{\partial y_{\alpha}} .
\end{equation} 

 For a Hamiltonian function $\mathcal{H}$ defined on the total space of the dual bundle $\mathcal{A}^{\ast }$ the Hamilton's equation is 
 \begin{equation}
     \frac{du}{dt}=X_{\mathcal{H}} (u)=\{u,\mathcal{H}\}_{\mathcal{A}^{\ast } }. 
 \end{equation}
 Accordingly, in finite dimensions, the Hamiltonian vector field $X_{\mathcal{H}}$ is computed to be 
\begin{equation}\label{H-Kesit}
X_{\mathcal{H}} = (\mathfrak{a}_{\mathcal{A}})_{\alpha }^{i}\frac{\partial \mathcal{H}}{\partial y_{\alpha}}\frac{\partial}{\partial x^{i}} - \left( \mathfrak{H}_{\alpha \beta}^{\gamma}y_{\gamma} \frac{\partial \mathcal{H}}{\partial y_{\beta}} + (\mathfrak{a}_{\mathcal{A}})_{\alpha }^{i}\frac{\partial \mathcal{H}}{\partial x^i} \right)\frac{\partial}{\partial y_{\alpha}}.
\end{equation}
The integral curves of the section given by this vector field satisfy Hamilton's equations in coordinates:
\begin{align}\label{ham-eq}
 \frac{d x^i}{dt} = (\mathfrak{a}_{\mathcal{A}})_{\alpha }^{i}\frac{\partial \mathcal{H}}{\partial y_{\alpha}},\qquad
 \frac{d y_{\alpha}}{dt} = - \mathfrak{H}_{\alpha \beta}^{\gamma}y_{\gamma} \frac{\partial \mathcal{H}}{\partial y_{\beta}} - (\mathfrak{a}_{\mathcal{A}})_{\alpha }^{i}\frac{\partial \mathcal{H}}{\partial x^i} .
\end{align}

\subsubsection*{Lagrangian Dynamics on Lie Algebroids.} 

Let $L=L(x^i, y^\alpha)$ be a Lagrangian function on the Lie algebroid $\mathcal{A}$. Referring the reader to \cite{martinez2001,WeinLag} for a coordinate-free study of the Lagrangian dynamics on Lie algebroids without referencing the Poisson structure on the dual bundle, we shall hereby consider the Legendre transformation (fiber derivative) as a map from the Lie algebroid to the dual bundle as 
\begin{equation}
\mathbb{F}L:\mathcal{A}\longrightarrow \mathcal{A}^\ast, \qquad (x^i, y^\alpha)\mapsto \big(x^i, y_\alpha=\frac{\partial L}{\partial y^\alpha}\big)
\end{equation}
for a more accessible way. Assuming that the Lagrangian is regular (non-degenerate) we argue that the Legendre transformation is a (at least local) diffeomorphism. Using the Legendre transformation, we pull the Poisson structure on the dual bundle $\mathcal{A}^\ast$ to the Lie algebroid, which yields a Poisson structure, denoted by $\{\bullet,\bullet\}_{\mathcal{A}}$, on the Lie algebroid. Define the energy function determined by $L$ as
\begin{equation}\label{energy-func}
E_L:\mathcal{A}\longrightarrow \mathbb{R}, \qquad  E_L(x^i,y^\alpha)=\frac{\partial L}{\partial y^\alpha}y^\alpha - L(x^i,y^\alpha) .
\end{equation}
 On the space $\mathcal{A}$, the dynamics is computed to be 
\begin{equation}
\frac{dx^i}{dt} = \{x^i, E_L\}_{\mathcal{A} }, \qquad \frac{dy^\alpha}{dt} = \{y^\alpha, E_L\}_{\mathcal{A} }. 
\end{equation}
Explicitly, the Euler-Lagrange equations on the Lie algebroid may be presented as follows, \cite{WeinLag}:
\begin{equation} \label{AEL}
\frac{d}{dt}\frac{\partial L}{\partial y^\alpha}= (\mathfrak{a}_{\mathcal{A}} )^i_\alpha \frac{\partial L}{\partial x^i} +  \mathfrak{H}^\gamma_{\alpha\beta}y^\beta\frac{\partial L}{\partial y^\gamma}, \qquad \frac{dx^i}{dt} = (\mathfrak{a}_{\mathcal{A}} )^i_\alpha y^\alpha.
\end{equation}

\subsubsection*{The Lie Algebra Framework.} 

In this paragraph, we shall examine the dynamical equations in the level of Lie algebras. To this end we let $M=\{e\}$, and accordingly the achor map to be trivial. Then we arrive at a Lie algebra $\mathfrak{g}$ as the space of sections. The dual space of this space is $\mathfrak{g}^*$. The Poisson bracket in \eqref{Poisson on dual A} reduces to so-called Lie-Poisson bracket. For two functions $f=f(u)$ and $h=h(u)$ defined on the dual space  $\mathfrak{g}^*$, the Lie-Poisson bracket is defined to be 
\begin{equation}
  \{f,h\} (u) = \Big\langle u   ,  \big[  \frac{\delta f} { \delta u}, 
\frac{\delta h} { \delta u} \big]
\Big\rangle  = \langle u   ,  ad_{ \delta f / \delta u} \frac{\delta h} { \delta u} 
\rangle,
\end{equation}
where the pairing is the natural one between $\mathfrak{g}^*$ and $\mathfrak{g}$ whereas the bracket inside pairing is the Lie algebra bracket on $\mathfrak{g}$. Here, $ad$ refers to the adjoint action of the Lie algebra on itself. We remark that  $\delta f/\delta u$ is the Fr\'echet derivative for the infinite dimensional case and partial derivative for the finite dimensional case. We also assume the reflexivity condition.   
In this reduced picture, the Hamiltonian vector field on $\mathfrak{g}^* $ generated by a Hamiltonian function $h$ is then computed to be
\begin{equation}\label{H-H-G-}
X_{h}(u)= - ad^*_{ \delta h / \delta u}u ,
\end{equation}
where $ad^*$ is the linear algebraic dual of the adjoint action $ad$. 
In the Lie-Poisson picture, the dynamics of an observable $f$, governed by a Hamiltonian function $h$, is computed to be
\begin{equation}\label{aaa} 
\dot{f}=\{f,h\} ( u  )
 =  - \Big\langle
ad_{\delta h/
\delta u }^{\ast }u  ,\frac{\delta f}{\delta u  }
\Big\rangle. 
\end{equation}
This bracket is exactly the one in \eqref{Poisson on dual A} where the Lie algebroid is chosen to be a Lie algebra. 
Then we have the Lie-Poisson equation on $\mathfrak{g}^*$ given by \index{Lie-Poisson equation }
\begin{equation}\label{LJ}
 \frac{d u}{dt}= - ad^*_{ \delta h / \delta u}u .
\end{equation} 
For finite dimensions, in coordinates $(y_\alpha)$ on $\mathfrak{g}^*$, the Lie-Poisson bracket is computed to be 
\begin{equation}
\{ y_{\alpha}, y_{\beta} \}_{\mathfrak{g}^*}  = \mathfrak{H}_{\alpha \beta}^{\gamma} y_{\gamma},
\end{equation} 
where in this case, $\mathfrak{H}_{\alpha \beta}^{\gamma}$ refers to the structure constants of the Lie algebra. Referring to the local Hamiltonian dynamics \eqref{ham-eq} on the Lie algebroid setting, we write the Lie-Poisson equations on the dual of a Lie algebroid as 
\begin{align}\label{ham-eq++} 
 \frac{d y_{\alpha}}{dt} = - \mathfrak{H}_{\alpha \beta}^{\gamma}y_{\gamma} \frac{\partial h}{\partial y_{\beta}}. 
\end{align}
From the Lagrangian dynamics point of view, we choose a Lagrangian function $l=l(X)$ on $\mathfrak{g}$. This reduces the Euler-Lagrange equations \eqref{AEL} on the Lie algebroid setting to the Euler-Poincar\'{e} equation
\begin{equation}
\frac{d}{dt}\frac{\delta l}{\delta X}= ad_{X}^{\ast}\frac{\delta l}{\delta X}.
\label{EP-eqn}
\end{equation}
on the Lie algebra. In finite dimensions, we choose $(y^\alpha)$ as coordinates on $\mathfrak{g}$, to obtain  
\begin{equation} \label{AEL++}
\frac{d}{dt}\frac{\partial l}{\partial y^\alpha}=   \mathfrak{H}^\gamma_{\alpha\beta}y^\beta\frac{\partial l}{\partial y^\gamma}.
\end{equation}

\subsection{Reversible Dynamics on BDCP Lie Algebroids} \label{sec-dyn-bicocycle}

This section is reserved for both the Lagrangian dynamics on bicocycle double cross product Lie algebroids, and the Hamiltonian dynamics on the dual space. This will give us both the collective motion of the dynamical systems as well as the decomposition. We start with the bicocycle double cross product  Poisson structure and the bicocycle double cross product Hamiltonian dynamics.  

\subsubsection*{BDCP Poisson Bracket.} 

We first recall the space of functions defined in \eqref{lolo} and \eqref{nono} and carry them to the dual space $(\mathcal{A}^{\ast}{\hspace{.1cm}}_{\zeta\hspace{-.1cm}} \bowtie _\psi \mathcal{B}^{\ast},\pi_{_\zeta\bowtie_\psi},M)$ of the unified product Lie algebroid structure $(\mathcal{A}_\zeta \bowtie_\psi \mathcal{B},\tau_{_\zeta\bowtie_\psi},M,\mathfrak{a}_{_\zeta \bowtie_\psi},[\bullet,\bullet]_{_\zeta \bowtie_\psi})$. The total space is topologically equivalent to $\mathcal{A}^{\ast} \times  \mathcal{B}^{\ast}$ but we denote it by $(\mathcal{A}^{\ast} _\zeta \bowtie_\psi \mathcal{B}^{\ast}, \pi_{_\zeta\bowtie_\psi},M)$ in order to highlight its product character.
The space of functions on the dual space $\mathcal{A}^{\ast} _\zeta \bowtie_\psi \mathcal{B}^{\ast}$ is spanned by three classes of functions. The first class of functions are 
\begin{equation}
\widehat{f}: \mathcal{A}^{\ast} {\hspace{.1cm}}_{\zeta\hspace{-.1cm}} \bowtie _\psi \mathcal{B}^{\ast} \longrightarrow \mathbb{R}, \qquad \widehat{f}(u) = f \circ \pi_{_\zeta \bowtie_\psi}(u).
\end{equation}
These functions are constants on the fibers. The second and third classes of functions come from the section of the individual spaces. Considering the sections $X$ in $\Gamma(\mathcal{A})$ and $Y$ in $\Gamma(\mathcal{B})$, we determine the linear independent sections of $\mathcal{A}_\zeta \bowtie_\psi \mathcal{B}$. We denote these sections as $X$ and $Y$ respectively, in order not to cause any notational inflation. That is, we take $X = (X,0)$ and $Y=(0,Y)$. They, then, give rise to functions on the total space as 
\begin{equation}
\widehat{Z}: \mathcal{A}^{\ast} {\hspace{.1cm}}_{\zeta\hspace{-.1cm}} \bowtie _\psi \mathcal{B}^{\ast} \longrightarrow \mathbb{R}, \qquad \widehat{Z}(w) = \langle w, Z \circ \pi_{ _\zeta \bowtie_\psi}(w) \rangle ,
\end{equation}
where $Z$ is either $X$ or $Y$. These functions are linear on the fibers. Therefore, in the light of the Poisson bracket in \eqref{Poisson on dual A}, we present the Poisson bracket on the dual space $\mathcal{A}^{\ast} {\hspace{.1cm}}_{\zeta\hspace{-.1cm}} \bowtie _\psi \mathcal{B}^{\ast}$ of the bicocycle double cross product Lie algebroid $\mathcal{A}{\hspace{.1cm}}_{\zeta\hspace{-.1cm}} \bowtie _\psi \mathcal{B}$ as
\begin{equation} \label{UP-Algoid-Poisson}
\begin{split}
    \{\widehat{f},\widehat{g}\}_{\mathcal{A}^{\ast} _\zeta \bowtie_\psi \mathcal{B}^{\ast}} & = 0,\\
    \{\widehat{f},\widehat{X}\}_{\mathcal{A}^{\ast} _\zeta \bowtie_\psi \mathcal{B}^{\ast}} & = \widehat{\mathcal{L}_{\mathfrak{a}_{_\zeta \bowtie_\psi} (X)} f},\\
    \{\widehat{f},\widehat{Y}\}_{\mathcal{A}^{\ast} _\zeta \bowtie_\psi \mathcal{B}^{\ast}} & = \widehat{\mathcal{L}_{\mathfrak{a}_{_\zeta \bowtie_\psi} (Y)} f}, \\
    \{\widehat{X}_1,\widehat{X}_2\}_{\mathcal{A}^{\ast} _\zeta \bowtie_\psi \mathcal{B}^{\ast}} & = \widehat{\phi(X_1,X_2)} + \widehat{\zeta(X_1,X_2)} \\
    \{\widehat{Y}_1,\widehat{Y}_2\}_{\mathcal{A}^{\ast} _\zeta \bowtie_\psi \mathcal{B}^{\ast}} & = \widehat{\psi(Y_1,Y_2)} + \widehat{\theta(Y_1,Y_2)}\\
    \{\widehat{Y},\widehat{X}\}_{\mathcal{A}^{\ast} _\zeta \bowtie_\psi \mathcal{B}^{\ast}} & = \widehat{\rho(Y,X)} - \widehat{\sigma(X,Y)},
\end{split}
\end{equation}
where one can consider the bracket on the fourth line as the result of the bracket of two elements $\widehat{X}_1$ and $\widehat{X}_2$ taking values both in $\mathcal{A}$ and $\mathcal{B}$, whereas the terms on the right side of the fifth line as the result of the bracket of two elements $\widehat{Y}_1$ and $\widehat{Y}_2$ taking values both in $\mathcal{A}$ and $\mathcal{B}$. The sixth line is the manifestation of the (weak) representations.

\subsubsection*{The Local Picture.} 

Let us now, for a brief moment, turn our attention to the general case of the dual space $\mathcal{A}^* {\hspace{.1cm}}_{\zeta\hspace{-.1cm}} \bowtie _\psi \mathcal{B}^*$ 
 of the bicocycle double cross product Lie algebroid $\mathcal{A} {\hspace{.1cm}}_{\zeta\hspace{-.1cm}} \bowtie _\psi \mathcal{B}$. We shall then examine the Poisson structure listed above. Let $(x^{i})$, $(x^{i},y^{\alpha })$, and $(x^{i},y_{\alpha })$ represents local coordinates on the base manifold $M$, on the Lie algebroid $\mathcal{A}$, and the dual space $\mathcal{A}^{\ast}$ respectively. 
On the vector bundle $\mathcal{B}$ we choose the local coordinates $(x^i,y^a)$, and on the dual space $\mathcal{B}^*$ we take $(x^i,y_a)$. As a result, on the space $\mathcal{A} {\hspace{.1cm}}_{\zeta\hspace{-.1cm}} \bowtie _\psi \mathcal{B} $  we have the induced local coordinates 
 $(x^i,y^\alpha,y^a)$, and on the  dual bundle $\mathcal{A}^*{\hspace{.1cm}}_{\zeta\hspace{-.1cm}} \bowtie _\psi \mathcal{B}^*$ we have $(x^i,y_\alpha,y_a)$. 
In this case, referring to the local characterizations \eqref{monon6},  the bicocycle double cross product Poisson bracket \eqref{Poisson on dual A} takes the form:
\begin{equation} \label{local Poisson on dual A}
\begin{split}
\{ x^{i}, x^{j} \}_{\mathcal{A}^{\ast} _\zeta \bowtie_\psi \mathcal{B}^{\ast}}  & = 0, \\ 
\{x^{i}, y_{\alpha} \}_{\mathcal{A}^{\ast} _\zeta \bowtie_\psi \mathcal{B}^{\ast}} &  = (\mathfrak{a}_{\mathcal{A}})_{\alpha}^{i}, \\ 
\{x^{i}, y_{a} \}_{\mathcal{A}^{\ast} _\zeta \bowtie_\psi \mathcal{B}^{\ast}} &  = (\mathfrak{a}_{\mathcal{B}})_{a}^{i},\\ 
\{ y_{\alpha}, y_{\beta} \}_{\mathcal{A}^{\ast} _\zeta \bowtie_\psi \mathcal{B}^{\ast}}  &= \mathfrak{H}_{\alpha \beta}^{\gamma} y_{\gamma} + \mathfrak{Z}_{\alpha \beta}^{a} y_{a}, 
\\ 
\{ y_{a}, y_{b} \}_{\mathcal{A}^{\ast} _\zeta \bowtie_\psi \mathcal{B}^{\ast}}  &= \mathfrak{P}_{a b}^{\gamma} y_{\gamma} + \mathfrak{T}_{a b}^{d} y_{d},
\\ 
\{ y_{a}, y_{\beta} \}_{\mathcal{A}^{\ast} _\zeta \bowtie_\psi \mathcal{B}^{\ast}} &= \mathfrak{R}_{a \beta}^{\gamma} y_{\gamma} + \mathfrak{S}_{a \beta}^{d} y_{d}.
\end{split}
\end{equation}
See that $\mathfrak{H}_{\alpha \beta}^{\gamma}$ is standing for the $\phi$ mapping, $\mathfrak{Z}_{\alpha \beta}^{a}$ is for the $\zeta$ mapping, $\mathfrak{P}_{a b}^{\gamma}$ is due to $\psi$, and $\mathfrak{T}_{a b}^{d}$ represents $\theta$.  

  Let us examine the Poisson bracket on the
 particular instances of the bicocycle double cross product construction. 
In case the $\zeta$ term is zero, we arrive at the unified product Poisson structure which we denote by $\{\bullet,\bullet\}_{\mathcal{A}^{\ast}  \bowtie_\psi \mathcal{B}^{\ast}}$. If both $\zeta$ and $\psi$ are zero, then we arrive at the double cross product Poisson structure $\{\bullet,\bullet\}_{\mathcal{A}^{\ast}  \bowtie  \mathcal{B}^{\ast}}$. If $\zeta$, $\theta$ and $\rho$ are all zero, then we have cocycle extension Poisson bracket $\{\bullet,\bullet\}_{\mathcal{A}^{\ast}  \ltimes _\psi\mathcal{B}^{\ast}}$. The case of the dual space of the double cross product Lie algebra has been studied in \cite{esen2016hamiltonian,EsSu21}.

\subsubsection*{Hamiltonian Dynamics on the Dual of BDCP Lie Algebroids.} 

Given a Hamiltonian function $\mathcal{H}=\mathcal{H}(x^i,y_\alpha,y_a)$ on the dual bundle $ \mathcal{A}^{\ast} {\hspace{.1cm}}_{\zeta\hspace{-.1cm}} \bowtie _\psi \mathcal{B}^{\ast}$, let us recall the Hamiltonian equation
\begin{equation}
\frac{du}{dt}=\{u,\mathcal{H}\}_{\mathcal{A}^{\ast} _\zeta \bowtie_\psi \mathcal{B}^{\ast}} 
\end{equation}
where the Poisson bracket on the right hand side is the one defined in \eqref{UP-Algoid-Poisson}. In local coordinates $(x^i,y_\alpha,y_a)$ on $ \mathcal{A}^{\ast} {\hspace{.1cm}}_{\zeta\hspace{-.1cm}} \bowtie _\psi \mathcal{B}^{\ast}$ and in view of the local Poisson brackets in \eqref{local Poisson on dual A}, we compute the Hamiltonian vector field as
\begin{equation}\label{H-Kesit+}
\begin{split}
X_{\mathcal{H}} & = \Big( (a_{\mathcal{B}})_{a}^{i}\frac{\partial \mathcal{H}}{\partial y_{a}}+(\mathfrak{a}_{\mathcal{A}})_{\alpha}^{i}\frac{\partial \mathcal{H}}{\partial y_{\alpha}}\Big) \frac{\partial}{\partial x^{i}} \\ 
& - \left( \mathfrak{Z}_{\alpha \beta}^{d}y_{d} \frac{\partial \mathcal{H}}{\partial y_{\beta}}+\mathfrak{H}_{\alpha \beta}^{\gamma}y_{\gamma} \frac{\partial \mathcal{H}}{\partial y_{\beta}}-\mathfrak{S}_{b \alpha}^{d}y_{d} \frac{\partial \mathcal{H}}{\partial y_{b}}-
\mathfrak{R}_{b \alpha}^{\gamma}y_{\gamma} \frac{\partial \mathcal{H}}{\partial y_{b}} + (\mathfrak{a}_{\mathcal{A}})_{\alpha }^{i}\frac{\partial \mathcal{H}}{\partial x^i} \right)\frac{\partial}{\partial y_{\alpha}}
\\ 
& - \left( 
\mathfrak{S}_{a \beta}^{d}y_{d} \frac{\partial \mathcal{H}}{\partial y_{\beta}}
+
\mathfrak{T}_{a b}^{d}y_{d} \frac{\partial \mathcal{H}}{\partial y_{b}} +
\mathfrak{P}_{a b}^{\gamma}y_{\gamma} \frac{\partial \mathcal{H}}{\partial y_{b}} +
\mathfrak{R}_{a \beta}^{\gamma}y_{\gamma} \frac{\partial \mathcal{H}}{\partial y_{\beta}}
+ (\mathfrak{a}_{\mathcal{A}})_{a}^{i}\frac{\partial \mathcal{H}}{\partial x^i} \right)\frac{\partial}{\partial y_{a}}.
\end{split}
\end{equation}
In accordance with this, what we call the \emph{bicocycle double cross product Hamilton's equations} are computed to be 
\begin{equation} \label{BDCPHamEqLoc}
\begin{split}
&\frac{d x^i}{dt} = (\mathfrak{a}_{\mathcal{B}})_{a}^{i}\frac{\partial \mathcal{H}}{\partial y_{a}}+(\mathfrak{a}_{\mathcal{A}})_{\alpha}^{i}\frac{\partial \mathcal{H}}{\partial y_{\alpha}},\\
& \frac{d y_{\alpha}}{dt} = -(\mathfrak{a}_{\mathcal{A}})_{\alpha }^{i}\frac{\partial \mathcal{H}}{\partial x^i}-   \mathfrak{H}_{\alpha \beta}^{\gamma}y_{\gamma} \frac{\partial \mathcal{H}}{\partial y_{\beta}}  +\underbrace{\mathfrak{S}_{b \alpha}^{d}y_{d} \frac{\partial \mathcal{H}}{\partial y_{b}}}_{\text{\textit{rep.} of }\mathcal{A} \text{ on } \mathcal{B}} -\underbrace{
\mathfrak{Z}_{\alpha \beta}^{d}y_{d} \frac{\partial \mathcal{H}}{\partial y_{\beta}}}_{\text{\textit{cocycle} term } \zeta } + \underbrace{
\mathfrak{R}_{b \alpha}^{\gamma}y_{\gamma} \frac{\partial \mathcal{H}}{\partial y_{b}}}_{\text{\textit{rep.}}  ~ of ~ \mathcal{B} ~  on ~  \mathcal{A}}  ,
\\
&
\frac{d y_{a}}{dt} =- (a_{\mathcal{B}})_{a}^{i}\frac{\partial \mathcal{H}}{\partial x^i} 
-
\mathfrak{T}_{a b}^{d}y_{d} \frac{\partial \mathcal{H}}{\partial y_{b}} - \underbrace{
\mathfrak{S}_{a \beta}^{d}y_{d} \frac{\partial \mathcal{H}}{\partial y_{\beta}}}_{\text{\textit{rep.} of }\mathcal{A} \text{ on } \mathcal{B}}
-\underbrace{
\mathfrak{P}_{a b}^{\gamma}y_{\gamma} \frac{\partial \mathcal{H}}{\partial y_{b}}}_{\text{\textit{cocycle} term } \psi }
- \underbrace{
\mathfrak{R}_{a \beta}^{\gamma}y_{\gamma} \frac{\partial \mathcal{H}}{\partial y_{\beta}}}_{\text{\textit{rep.}}  ~ of ~ \mathcal{B} ~  on ~  \mathcal{A}} 
.
\end{split}
\end{equation}

\subsubsection*{The (Dual) Lie Algebra Framework.} 

Having formulated them for bicocycle double cross product Lie algebroids, we can now present the ingredients of the Poisson and the Hamiltonian dynamics for bicocycle double cross product Lie algebras without much difficulty. In order to achieve this we assume the base manifold to be a singleton, and the anchor maps to be trivial. In this case, one has two vector spaces $\mathfrak{g}$ and $\mathfrak{h}$ and their bicocycle double cross product Lie algebra $\mathfrak{g}{\hspace{.1cm}}_{\zeta\hspace{-.1cm}} \bowtie _\psi \mathfrak{h}$ equipped with the Lie bracket $[\bullet,\bullet]_{_\zeta \bowtie_\psi}$. Correspondingly, the dual space $\mathfrak{g}^{\ast}{\hspace{.1cm}}_{\zeta\hspace{-.1cm}} \bowtie _\psi \mathfrak{h}^{\ast}$ is endowed with the bicocycle double cross product Lie-Poisson bracket, which in local coordinates is given by 
\begin{equation} \label{local Poisson on dual A ++}
\begin{split}
\{ y_{\alpha}, y_{\beta} \}_{\mathfrak{g}^{\ast} _\zeta \bowtie_\psi \mathfrak{h}^{\ast}}  &= \mathfrak{H}_{\alpha \beta}^{\gamma} y_{\gamma} + \mathfrak{Z}_{\alpha \beta}^{a} y_{a}, 
\\ 
\{ y_{a}, y_{b} \}_{\mathfrak{g}^{\ast} _\zeta \bowtie_\psi \mathfrak{h}^{\ast}}  &= \mathfrak{P}_{a b}^{\gamma} y_{\gamma} + \mathfrak{T}_{a b}^{d} y_{d},
\\ 
\{ y_{a}, y_{\beta} \}_{\mathfrak{g}^{\ast} _\zeta \bowtie_\psi \mathfrak{h}^{\ast}} &= \mathfrak{R}_{a \beta}^{\gamma} y_{\gamma} + \mathfrak{S}_{a \beta}^{d} y_{d},
\end{split}
\end{equation}
where the coefficients are those that determine the representations / twisted cocycles just as in the Lie algebroid setting. Given Hamiltonian function $h=h(y_{\alpha},y_a)$ on the dual space, the bicocycle double cross product Lie-Poisson equations are then computed to be 
\begin{equation} \label{BDCPHamEqLoc++}
\begin{split}
& \frac{d y_{\alpha}}{dt} = -   \mathfrak{H}_{\alpha \beta}^{\gamma}y_{\gamma} \frac{\partial h}{\partial y_{\beta}}  +\underbrace{\mathfrak{S}_{b \alpha}^{d}y_{d} \frac{\partial h}{\partial y_{b}}}_{\text{\textit{rep.} of }\mathfrak{g}  \text{ on } \mathfrak{h} } -\underbrace{
\mathfrak{Z}_{\alpha \beta}^{d}y_{d} \frac{\partial h}{\partial y_{\beta}}}_{\text{\textit{cocycle} term } \zeta } + \underbrace{
\mathfrak{R}_{b \alpha}^{\gamma}y_{\gamma} \frac{\partial h}{\partial y_{b}}}_{\text{\textit{rep.}}  ~ of ~ \mathfrak{h} ~  on ~  \mathfrak{g} }  ,
\\
&
\frac{d y_{a}}{dt} =
-
\mathfrak{T}_{a b}^{d}y_{d} \frac{\partial h}{\partial y_{b}} - \underbrace{
\mathfrak{S}_{a \beta}^{d}y_{d} \frac{\partial h}{\partial y_{\beta}}}_{\text{\textit{rep.} of }\mathfrak{g}  \text{ on } \mathfrak{h}}
-\underbrace{
\mathfrak{P}_{a b}^{\gamma}y_{\gamma} \frac{\partial h}{\partial y_{b}}}_{\text{\textit{cocycle} term } \psi }
- \underbrace{
\mathfrak{R}_{a \beta}^{\gamma}y_{\gamma} \frac{\partial h}{\partial y_{\beta}}}_{\text{\textit{rep.}}  ~ of ~ \mathfrak{h} ~  on ~  \mathfrak{g} } 
.
\end{split}
\end{equation}

\subsubsection*{Lagrangian  Dynamics on the BDCP Lie Algebroids.}

Let us fix a local coordinate system on the Lie algebroid $\mathcal{A}{\hspace{.1cm}}_{\zeta\hspace{-.1cm}} \bowtie _\psi \mathcal{B}$ as $(x^i, \bar{y}^{\kappa}) := (x^i, y^\alpha, y^a)$. Similar to the Euler-Lagrange equations on a Lie algebroid, the Euler-Lagrange equations generated by a Lagrangian function $L$ on the product Lie algebroid $\mathcal{A}{\hspace{.1cm}}_{\zeta\hspace{-.1cm}} \bowtie _\psi \mathcal{B}$ (in local coordinates) are two sets of equations. The first set of  equations can be directly interpreted as:
\begin{equation}\label{BDCP-EL1}
\frac{dx^i}{dt} =  (\mathfrak{a}_{\mathcal{A}})^i_\alpha \,y^\alpha+(\mathfrak{a}_{\mathcal{B}})^i_a y^a.
\end{equation}
The second set of equations can be expressed as:
\begin{equation} \label{BDCP-EL}
\begin{split}
&  \frac{d}{dt}\frac{\partial L}{\partial y^\beta} =(\mathfrak{a}_{\mathcal{A}})^i_\beta\frac{\partial L}{\partial x^i} + \mathfrak{H}^\alpha_{\beta\gamma}y^\gamma\frac{\partial L}{\partial y^\alpha} -
\underbrace{\mathfrak{R}^\alpha_{d \beta }y^d \frac{\partial L}{\partial y^\alpha}}_{\text{\textit{rep.}}  ~ of ~ \mathcal{B} ~  on ~  \mathcal{A}} -
\underbrace{\mathfrak{S}^a_{d\beta }y^d\frac{\partial L}{\partial y^a}}_{\text{\textit{rep.} of }\mathcal{A} \text{ on } \mathcal{B}}  \\
& \qquad \qquad\qquad\qquad +
\underbrace{\mathfrak{Z}^a_{\beta \gamma}y^\gamma\frac{\partial L}{\partial y^a}}_{\text{\textit{cocycle} term } \zeta },\\
&  \frac{d}{dt}\frac{\partial L}{\partial y^b} = (a_\mathcal{B})^i_b \frac{\partial L}{\partial x^i} +\mathfrak{T}^a_{bd}y^d\frac{\partial L}{\partial y^a} + 
\underbrace{\mathfrak{S}^a_{b\gamma}y^\gamma\frac{\partial L}{\partial y^a}}_{\text{\textit{rep.} of }\mathcal{A} \text{ on } \mathcal{B}}+ \underbrace{\mathfrak{R}^\alpha_{b\gamma}y^\gamma\frac{\partial L}{\partial y^\alpha}}_{\text{\textit{rep.}}  ~ of ~ \mathcal{B} ~  on ~  \mathcal{A}}  \\
& \qquad \qquad\qquad\qquad + \underbrace{\mathfrak{P}^\alpha_{bd}y^d\frac{\partial L}{\partial y^\alpha}}_{\text{\textit{cocycle} term } \psi } .
\end{split}
\end{equation}

\subsubsection*{The Lie Algebra Framework.} 

Let us examine the Euler-Poincar\'{e} dynamics on the bicocycle double cross product Lie algebra $\mathfrak{g}{\hspace{.1cm}}_{\zeta\hspace{-.1cm}} \bowtie _\psi \mathfrak{h}$. To this end, we set the terms responsible for the anchor maps in \eqref{BDCP-EL} to be trivial. Then, given a Lagrangian function $l=l(y^\alpha,y^a)$, in local coordinates $(y^\alpha,y^a)$, the bicocycle double cross product Euler-Poincar\'{e} equations take the form 
\begin{equation} \label{BDCP-EP}
\begin{split}
&  \frac{d}{dt}\frac{\partial l}{\partial y^\beta} = \mathfrak{H}^\alpha_{\beta\gamma}y^\gamma\frac{\partial l}{\partial y^\alpha} -
\underbrace{\mathfrak{R}^\alpha_{d \beta }y^d \frac{\partial l}{\partial y^\alpha}}_{\text{\textit{rep.}}  ~ of ~ \mathfrak{h}  ~  on ~  \mathfrak{g} } -
\underbrace{\mathfrak{S}^a_{d\beta }y^d\frac{\partial l}{\partial y^a}}_{\text{\textit{rep.} of }\mathfrak{g}  \text{ on } \mathfrak{h} }   +
\underbrace{\mathfrak{Z}^a_{\beta \gamma}y^\gamma\frac{\partial l}{\partial y^a}}_{\text{\textit{cocycle} term } \zeta },\\
&  \frac{d}{dt}\frac{\partial l}{\partial y^b} =  \mathfrak{T}^a_{bd}y^d\frac{\partial l}{\partial y^a} + 
\underbrace{\mathfrak{S}^a_{b\gamma}y^\gamma\frac{\partial l}{\partial y^a}}_{\text{\textit{rep.} of }\mathfrak{g}  \text{ on } \mathfrak{h} }+ \underbrace{\mathfrak{R}^\alpha_{b\gamma}y^\gamma\frac{\partial l}{\partial y^\alpha}}_{\text{\textit{rep.}}  ~ of ~ \mathfrak{h}  ~  on ~  \mathfrak{g} }  + \underbrace{\mathfrak{P}^\alpha_{bd}y^d\frac{\partial l}{\partial y^\alpha}}_{\text{\textit{cocycle} term } \psi } .
\end{split}
\end{equation}

\subsection{Irreversible Dynamics via Real Line Extension}\label{cont-sec}

Hamiltonian dynamics is defined on Poisson (more specifically symplectic) manifolds that offer a reversible motion, where the Hamiltonian function is preserved. A generalization of the classical Hamiltonian dynamics that allows the irreversible dynamics is achieved on contact manifolds. In the present subsection, we shall recall this formalism.

\subsubsection*{Contact Manifolds.} 

Let $N$ be an odd-dimensional manifold, say, of dimension $(2n + 1)$. A contact structure on $N$ is a maximally non-integrable smooth distribution of codimension one, defined locally by a one-form $\eta$ with the non-integrability condition 
\begin{equation}
d\eta^n \wedge \eta \neq 0.
\end{equation}
Such a one-form $\eta$ is called a (local) contact form \cite{Arn,Liber87}. We consider the existence of a global contact one-form and denote a contact manifold by $(N, \eta)$. Given a contact one-form $\eta$, the vector field $\mathcal{R}$ satisfying 
\begin{equation}
\iota_{\mathcal{R}}\eta =1,\qquad \iota_{\mathcal{R}}d\eta =0
\end{equation}
is unique, and is called the Reeb vector field. 
For a contact manifold $(N, \eta)$, there is a musical isomorphism 
\begin{equation}\label{flat-map-cont}
\flat:\mathfrak{X}(N)\longrightarrow \Gamma^1(N),\qquad \xi\mapsto \iota_\xi d\eta+\langle \eta, \xi  \rangle \eta.
\end{equation} 
We denote the inverse of \eqref{flat-map-cont} by $\sharp$.  
Referring to this, we define a bivector field $\Lambda$ on  $N$ as
\begin{equation}\label{Lambda}
\Lambda(\alpha,\beta)= d\eta(\sharp\alpha, \sharp \beta). 
\end{equation}
The pair $(\Lambda,\mathcal{R})$ gives a Jacobi structure \cite{Lichnerowicz-Jacobi,Marle-Jacobi}, since they satisfy \eqref{ident-Jac}. This realization allows us to define a Jacobi (contact) bracket 
\begin{equation}\label{cont-bracket-L+} 
\{F,\mathcal{H}\}^c =\Lambda(dF,d\mathcal{H}) +F\mathcal{R}(\mathcal{H}) - \mathcal{H}\mathcal{R}(F).
\end{equation}
We do note that the above bracket satisfies the Jacobi identity, but not the Leibniz identity.

\subsubsection*{Dissipation on Extended Manifolds.} 

Let $(N,\eta)$ be a a contact manifold. We define the contact Hamiltonian vector field $\xi _{\mathcal{H}}$ as \cite{Br17,BrCrTa17,de2021hamilton,de2019contact,esen2021contact}
\begin{equation}\label{Ham-def} 
\iota_{\xi _{\mathcal{H}}}d\eta=d\mathcal{H}-\mathcal{R}(\mathcal{H})\eta,\qquad \iota_{\xi _{\mathcal{H}}}\eta=-\mathcal{H}.
\end{equation} 
The relationship between the Hamiltonian vector field and the contact bracket is given by the equation
\begin{equation}\label{bra-X-J}
\xi_\mathcal{H}(F)= \{F,\mathcal{H}\}^{J}- F \mathcal{R}(\mathcal{H}).
\end{equation} 
Given Hamiltonian function $\mathcal{H}$, the contact Hamiltonian vector field is a strict contact vector field if and only if $\mathcal{R}(\mathcal{H})=0$. Accordingly, one has strict contact Hamiltonian vector field as
\begin{equation}\label{X-qua-ham}
  \iota_{\xi_{\mathcal{H}}}d\eta =d\mathcal{H},\qquad\iota_{\xi_{\mathcal{H}}}\eta =-\mathcal{H}  .
\end{equation}
In order to observe the dissipative behaviour of this dynamics we exhibit two calculations:
\begin{equation}
{\mathcal{L}}_{\xi_{\mathcal{H}}}  \mathcal{H} = - \mathcal{R}(\mathcal{H}) \mathcal{H},\qquad {\mathcal{L}}_{\xi_{\mathcal{H}}}  (\eta \wedge (d \eta)^n) = - (n+1) \, \mathcal{R}^{c} (\mathcal{H}) \eta \wedge (d\eta)^n,
\end{equation}
where the former one is showing that the Hamiltonian function is not preserved under the Hamiltonian function, whereas the latter gives the dissipation of the volume form $\eta \wedge (d \eta)^n$ (given for $(2n+1)$-dimenisonal contact manifold). Accordingly, the divergence of the contact Hamiltonian vector field is computed to be 
\begin{equation} \label{div-X-H}
\mathrm{div}(\xi_{\mathcal{H}})= -  (n+1)  \mathcal{R}(\mathcal{H}).
\end{equation}

Let now $Q$ be a manifold, and let $T^*Q$ be its cotangent bundle equipped with the canonical symplectic two-form $\Omega_Q$ defined to be the exterior derivative of the canonical one-form $\theta_Q$ with the negative sign. In Darboux coordinates $(x^i,p_i)$ on the cotangent bundle $T^*Q$, the canonical one-form and the canonical symplectic two-form can thus be presented as
\begin{equation}
\Theta_Q=p_i dx^i,\qquad \Omega_Q=-d\Theta_Q=dx^i\wedge dp_i
\end{equation}
respectively. 
A generic example of a contact manifold is the extended symplectic manifold $N=T^*Q\times \mathbb{R}$. It is to write the contact one-form referring to the canonical one-form $\Theta_Q$ and the differential one-form $dz$ defined on $\mathbb{R}$. More precisely, by the abuse of the notation (we denote pullbacks and pushforwards of the tensor fields by the same notation not to cause notation inflation), we have
\begin{equation}
    \eta=dz-\Theta_Q.
\end{equation}
This construction will be our main interest in the sequel, as it involves the addition of the dissipative terms to the dynamics on Lie algebroids. 

In the Darboux coordinates $(x^i,p_i,z)$ on $N$, the contact one-form and the Reeb field are computed to be
\begin{equation}\label{eta-M}
    \eta=dz-p_idx^i,\qquad \mathcal{R}=\frac{\partial}{\partial z}.
\end{equation}
Furthermore, in the Darboux coordinates, the bivector $\Lambda_{T^*Q\times \mathbb{R}}$ is computed to be 
\begin{equation}\label{AP-Bra}
\Lambda_{T^*Q\times \mathbb{R}}= \frac{\partial  }{\partial x^{i}}\wedge \frac{\partial  }{\partial p_{i}} + p_i \frac{\partial  }{\partial z} \wedge \frac{\partial  }{\partial p_i}. 
\end{equation} 
Here, the first term is due to the canonical Poisson bivector field $\Lambda_{T^*Q} $ on $T^*Q$ obtained from the symplectic two-form $\Omega_Q$. By defining the Liouville vector field 
\begin{equation}
\Delta^{*} (x^{i},p_{i},z)=p_{i}\frac{\partial}{\partial p_{i}}
\end{equation}
 on the extended cotangent bundle, we can write the contact bivector field as 
	\begin{equation}\label{contactization}
	\Lambda_{T^*Q\times \mathbb{R}}= \Lambda_{T^*Q} +\mathcal{R}\wedge  \Delta^{*}.  
	\end{equation}
In the local picture, the contact bracket \eqref{cont-bracket-L+} is written as
\begin{equation}\label{Lag-Bra}
\{F,\mathcal{H}\}^c  = \frac{\partial F}{\partial x^{i}}\frac{\partial \mathcal{H}}{\partial p_{i}} -
\frac{\partial F}{\partial p_{i}}\frac{\partial \mathcal{H}}{\partial x^{i}} + \left(F  - p_{i}\frac{\partial F}{\partial p_{i}} \right)\frac{\partial \mathcal{H}}{\partial z} -
\left(\mathcal{H}  - p_{i}\frac{\partial \mathcal{H}}{\partial p_{i}} \right)\frac{\partial F}{\partial z}.
\end{equation} 

Hence, the contact Hamiltonian vector field for a given Hamiltonian function $\mathcal{H}=\mathcal{H}
(x^i,p_i,z)$ becomes
\begin{equation}\label{con-dyn-X}
\xi_{\mathcal{H}}  =\frac{\partial \mathcal{H}}{\partial p_{i}}\frac{\partial}{\partial x^{i}}  - \left(\frac{\partial \mathcal{H}}{\partial x^{i}} + \frac{\partial \mathcal{H}}{\partial z} p_{i} \right)
\frac{\partial}{\partial p_{i}} + \left(p_{i}\frac{\partial \mathcal{H}}{\partial p_{i}} - \mathcal{H}\right)\frac{\partial}{\partial z}.
\end{equation}
We thus obtain the contact Hamilton's equations for $\mathcal{H}$ as
\begin{equation}\label{conham-X}
\frac{dx^i}{dt} = \frac{\partial \mathcal{H}}{\partial p_{i}}, 
\qquad
\frac{dp_i}{dt}= -\frac{\partial \mathcal{H}}{\partial x^{i}}- 
p_{i}\frac{\partial \mathcal{H}}{\partial z}, 
\qquad \frac{dz}{dt}= p_{i}\frac{\partial \mathcal{H}}{\partial p_{i}} - \mathcal{H}.
\end{equation}

\subsubsection*{Dissipative Lagrangian Dynamics.} 

In order to facilitate dissipation within the Lagrangian dynamics, on the tangent bundle $TQ$, we extend the tangent bundle by the real line. We thus arrive at the  extended tangent bundle  $T Q \times \mathbb{R}$, equipped with the induced coordinates $(x^i,\dot{x}^i,z)$. Hence, the dissipative dynamics is generated by a Lagrangian function $L=L(x^i,\dot{x}^i,z)$ on $TQ \times \mathbb{R}$. Here, we assume $(x^i,\dot{x}^i,z)$ as the coordinates on the extended tangent bundle.  
To derive the dynamical equations governed by such a Lagrangian function, the application of the Herglotz principle is necessary. This principle is defined by an action functional, as described in \cite{de2020review, Guenther, Herglotz, herg2}, where we refer the reader to \cite{de2019singular} for a more comprehensive account involving singular Lagrangians. In the present work, we shall stick to the regular Lagrangians. 

For a regular Lagrangian function $L$, the fiber derivative determines a diffeomorphism from the extended tangent bundle $T Q \times \mathbb{R}$ to the extended cotangent bundle $T^* Q \times \mathbb{R}$ as
\begin{equation}\label{Leg-Trf}
\mathbb{F}L: T Q \times \mathbb{R} \longrightarrow  T^* Q \times \mathbb{R}, \qquad (x^i,\dot{x}^i,z)\mapsto 
(x^i,\frac{\partial L}{\partial \dot{x}^i},z).
\end{equation}
We pull the contact one-form $\eta$ on $T^* Q \times \mathbb{R}$ and write it as $\eta^T$. This reads a contact manifold $(T Q \times \mathbb{R},\eta^T)$. The energy function is then defined to be
\begin{equation}
E_L(x^i,\dot{x}^i,z)= \dot{x}^i \frac{\partial L}{\partial \dot{x}^i} - L (x^i,\dot{x}^i,z),
\end{equation}
via which we determine the contact motion on the extended tangent bundle $TQ\times \mathbb{R}$. 
On the dual picture, the energy function corresponds to the Hamiltonian function
\begin{equation}
\mathcal{H}(x^i,p_i,z)=\dot{x}^ip_i-L(x^i,\dot{x}^i,z).
\end{equation} 
Hence, the contact dynamics \eqref{X-qua-ham} given in coordinate-free manner, assuming $\mathcal{H}:=E_L$ and $\eta:=\eta^T$ as the contact Hamiltonian function and the contact one-form, respectively, we arrive at the Herglotz equations  (also known as the generalized Euler--Lagrange equations) 
\begin{equation}\label{clagrangian4}
\frac{d}{dt} \frac{\partial L}{\partial \dot{x}^i} - \frac{\partial L}{\partial x^i} =
\frac{\partial L}{\partial \dot{x}^i} \frac{\partial L}{\partial z}.
\end{equation}
Removing the explicit dependence of $z$ reduces it to the classical Euler-Lagrange equations.

\subsection{Irreversible Dynamics on Lie Algebroids}\label{sec-irr-Lie}

Along the lines of  \cite{anahory2023reduction,anahory2024euler,Simoes2020b}, let $(\mathcal{A},\tau,M,\mathfrak{a}_{\mathcal{A}},[\bullet,\bullet]_{\mathcal{A}})$ be a Lie algebroid and $(\mathcal{A}^*,\pi,M)$ is its dual vector bundle equipped with the associated Poisson structure $\Lambda_{\mathcal{A}^*}$ defined globally in \eqref{Poisson on dual A}, which is locally given by \eqref{Poisson-A*}. Since the dissipative dynamics is formulated by the canonical Poisson structure of the real line $\mathbb{R}$ extension, we extend the Poisson manifold $\mathcal{A}^*$ by $\mathbb{R}$ to define the \emph{contactization} of the Poisson structure $\Lambda_{\mathcal{A}^*}$. We thus arrive at a Jacobi manifold structure on $\mathcal{A}^*\times \mathbb{R}$. Accordingly, setting the local coordinates $(x^i,y_\alpha,z)$ on the extended bundle $\mathcal{A}^*\times \mathbb{R}$, and referring to \eqref{contactization}, we define the bivector field 
\begin{equation} 
\Lambda_{\mathcal{A}^{*}\times \mathbb{R}} = \Lambda_{\mathcal{A}^{*}} + \mathcal{R} \wedge\Delta^{*}  = \frac{1}{2} \mathfrak{H}_{\alpha \beta}^{\gamma} y_{\gamma} \frac{\partial}{\partial y_{\alpha}} \wedge \frac{\partial}{\partial y_{\beta}} + \frac{\partial}{\partial x^{i}}\wedge \frac{\partial}{\partial y_{\alpha}} + y_{\alpha} \frac{\partial  }{\partial z} \wedge \frac{\partial  }{\partial y_{\alpha}}
\end{equation}
where the Reeb field and the Liouville vector field are  given by
 \begin{equation}
    \mathcal{R}=\frac{\partial}{ \partial z},\qquad \Delta^{*} = y_\alpha   \frac{\partial}{ \partial y_\alpha }, 
 \end{equation} 
 respectively. 
It then follows that the pair $(\Lambda_{\mathcal{A}^{*}\times \mathbb{R}},\mathcal{R})$ determines a Jacobi structure. As such, we can define a Jacobi bracket on $\mathcal{A}^{*}\times \mathbb{R}$ as   
\begin{equation}\label{cont-bracket-L} 
\{F,\mathcal{H}\}^j =\Lambda_{\mathcal{A}^{*}\times \mathbb{R}}(dF,dH) +F\mathcal{R}(\mathcal{H}) - \mathcal{H}\mathcal{R}(F).
\end{equation}
Assuming the local coordinates $(x^i,y_\alpha,z)$ on the extended space, the local realization of this bracket is computed as 
\begin{equation}\label{Pois-bra-A*xR}
	\begin{split}
& \{ x^{i}, x^{j} \}_{\mathcal{A}^{\ast}\times \mathbb{R}}   = 0, \\
		&  \{ z,x^{i}\}_{\mathcal{A}^{*}\times \mathbb{R}}=0,\\
		& \{z,y_\alpha\}_{\mathcal{A}^{*}\times \mathbb{R}}=y_\alpha, \\&  
\{x^{i}, y_{\alpha} \}_{\mathcal{A}^{\ast}\times \mathbb{R}}  = (\mathfrak{a}_{\mathcal{A}})_{\alpha}^{i},  \\&   
\{ y_{\alpha}, y_{\beta} \}_{\mathcal{A}^{\ast}\times \mathbb{R}}  = \mathfrak{H}_{\alpha \beta}^{\gamma} y_{\gamma},
	\end{split}	
\end{equation}
where $(\mathfrak{a}_{\mathcal{A}})_{\alpha}^{i}$ is the anchor map, whereas $\mathfrak{H}_{\alpha \beta}^{\gamma}$ are the structure functions of the Lie algebroid.

\subsubsection*{Hamiltonian Dynamics on Extended Lie Algebroids.}

Given a Hamiltonian function $\mathcal{H}=\mathcal{H}(x^i,y_\alpha,z)$ on the extended dual space $\mathcal{A}^{\ast}\times \mathbb{R}$, the Hamiltonian vector field  is defined to be
\begin{equation}\label{bra-X-J-}
\xi_\mathcal{H}  (F)= \{F,\mathcal{H}\}^{j}- F \mathcal{R}(\mathcal{H}).
\end{equation} 
So, in the local coordinates $(x^i,y_\alpha,z)$ on $\mathcal{A}^{\ast}\times \mathbb{R}$, the Hamiltonian vector field is computed to be  
\begin{equation}
	\xi_\mathcal{H}  =(\mathfrak{a}_{\mathcal{A}})_{\alpha}^{i}\frac{\partial \mathcal{H}}{\partial y_{\alpha}}\frac{\partial}{\partial x^{i}}-\Big( (\mathfrak{a}_{\mathcal{A}})_{\alpha}^{i}\frac{\partial \mathcal{H}}{\partial x^{i}}+ \mathfrak{H}_{\alpha \beta}^{\gamma} y_{\gamma}\frac{\partial \mathcal{H}}{\partial y_{\beta}}+ y_{\alpha}\frac{\partial \mathcal{H}}{\partial z}\Big)\frac{\partial}{\partial y_{\alpha}}+\left( y_{\alpha}\frac{\partial \mathcal{H}}{\partial y_{\alpha}}-\mathcal{H} \right)\frac{\partial}{\partial z}.
\end{equation}
Consequently, the dissipative Hamiltonian dynamics on the extended dual space $\mathcal{A}^{\ast}\times \mathbb{R}$, of the Lie algebroid $\mathcal{A}$, is given by 
\begin{equation}
\begin{split}
	\frac{dx^i}{dt} &=(\mathfrak{a}_{\mathcal{A}})_{\alpha}^{i}\frac{\partial \mathcal{H}}{\partial y_{\alpha}} , \\ \frac{dy_\alpha}{dt}&=-  (\mathfrak{a}_{\mathcal{A}})_{\alpha}^{i}\frac{\partial \mathcal{H}}{\partial x^{i}}- \mathfrak{H}_{\alpha \beta}^{\gamma} y_{\gamma}\frac{\partial \mathcal{H}}{\partial y_{\beta}}- y_{\alpha}\frac{\partial \mathcal{H}}{\partial z} ,\\ \frac{dz}{dt}&=  y_{\alpha}\frac{\partial \mathcal{H}}{\partial y_{\alpha}}-\mathcal{H}    .
 \end{split}
\end{equation}
A particular choice is $\mathcal{A}=TM$, and the anchor map to be the identity. In this case, $\mathcal{A}^*=T^*M$, the Hamiltonian vector field reduces to the contact Hamiltonian vector field \eqref{con-dyn-X}, and the dynamical equations turn out to be the contact Hamilton's equations \eqref{conham-X}. 

On the other extreme, if we choose the base manifold $M=\{e\}$ as a singleton, then the Lie algebroid $\mathcal{A}$ becomes a Lie algebra $\mathfrak{g}$, which will be investigated in greater detail below.

\subsubsection*{Lagrangian Dynamics on Extended Lie Algebroids.}

We shall now obtain the Lagrangian picture of the real line extension that controls the dissipative terms. Although, in general, it is possible to formulate the theory for any (even singular) Lagrangian function, we shall hereby confine ourselves with the regular Lagrangians.  

Let, now, $L=L(x^i,y^\alpha,z)$ be Lagrangian function on the extended Lie algebroid $\mathcal{A}\times \mathbb{R}$ with (local) coordinates $(x^i,y^\alpha,z)$. We then define the fiber derivative (or the Legendre transformation)  
\begin{equation} 
		\mathbb{F} L: \mathcal{A} \times \mathbb{R}  \longrightarrow \mathcal{A}^{*}\times \mathbb{R},\qquad 
		(x^i,y^\alpha,z)  \mapsto 	(x^i,\frac{\partial L}{\partial y^\alpha},z,
\end{equation}
which is at least a local diffeomorphism as a result of the regularity of the Lagrangian function. This, then, allows us to pull the Jacobi structure $(\Lambda_{\mathcal{A}^{*}\times \mathbb{R}},\mathcal{R})$ back to the extended Lie algebroid. We denote the Jacobi structure on $\mathcal{A} \times \mathbb{R}$ by the pair $(\Lambda_{\mathcal{A}\times \mathbb{R}},\mathcal{R}) $. This pair determines a Jacobi bracket $\{\bullet,\bullet\}_{\mathcal{A}\times \mathbb{R}}$ on the extended space. Defining the energy function 
\begin{equation}
    E_L(x^i,y^\alpha,z)= y^\alpha\frac{\partial L}{\partial y^\alpha} - L(x^i,y^\alpha,z)
\end{equation}
as the Hamiltonian function, we can determine the equations of motion as  in the form of Herglotz equations
	\begin{equation}\label{HerlAlg}
 \begin{split}
& \frac{d}{dt}\frac{\partial L}{\partial y^{\alpha}}-(\mathfrak{a}_{\mathcal{A}})_{\alpha}^{i}\frac{\partial L}{\partial x^{i}}=\mathfrak{H}_{\alpha \beta}^{\gamma}y^{\beta}\frac{\partial L}{\partial y^{\gamma}} +\frac{\partial L}{\partial z}\frac{\partial L}{\partial y^{\alpha}}, \\ & \frac{dx^i}{dt} = (\mathfrak{a}_{\mathcal{A}})_{\alpha}^{i}y^{\alpha} ,  \\ & \frac{dz}{dt}=L(x^i,y^\alpha,z). 
 \end{split}
	\end{equation} 
Note that in the case of the Lie algebroid $\mathcal{A}=TM$, the equations \eqref{HerlAlg} are just Herglotz equations \eqref{HerlAlg}
to the classical form in \eqref{clagrangian4}.

\subsubsection*{Lie Algebra Framework.} 

In the case of the base manifold is a point, the extended dual space happens to be $\mathfrak{g}^*\times \mathbb{R}$. Then assuming the coordinates $(y_\alpha , z)$, the Jacobi bracket may be presented as
\begin{equation}\label{Pois-bra-A*xR++}
	\begin{split}
 		 \{z,y_\alpha\}_{\mathfrak{g}^*\times \mathbb{R}}&=y_\alpha,   \\   
\{ y_{\alpha}, y_{\beta} \}_{\mathfrak{g}^*\times \mathbb{R}}  &= \mathfrak{H}_{\alpha \beta}^{\gamma} y_{\gamma}.
	\end{split}	
\end{equation}
Hence, given a Hamiltonian function $h=h(y_\alpha , z)$ on $\mathfrak{g}^*\times \mathbb{R}$, the Lie-Poisson-Herglotz equations are obtained as
\begin{equation}
\begin{split}
     \frac{dy_\alpha}{dt} &=-  \mathfrak{H}_{\alpha \beta}^{\gamma} y_{\gamma}\frac{\partial h}{\partial y_{\beta}} - y_{\alpha}\frac{\partial h}{\partial z} ,\\ \frac{dz}{dt}&=  y_{\alpha}\frac{\partial h}{\partial y_{\alpha}}-h   .
\end{split}
\end{equation}
As for a Lagrangian function $l=l(y^\alpha,z)$ on $\mathfrak{g}\times \mathbb{R}$, the equations \eqref{HerlAlg} reduces to the Euler-Poincar\'e-Herglotz equations
\begin{equation}\label{EPHerlAlg} 
\begin{split}
&\frac{d}{dt}\frac{\partial l}{\partial y^{\alpha}}+\mathfrak{H}_{\alpha \beta}^{\gamma}y^{\beta}\frac{\partial l}{\partial y^{\gamma}}=\frac{\partial l}{\partial z}\frac{\partial l}{\partial y^{\alpha}},  \\ &\frac{dz}{dt}=l(y^\alpha,z). 
\end{split}
	\end{equation}

\subsection{Irreversible  Dynamics on Bicocycle Double Cross Product Lie Algebroids} \label{noveldyn}

The present subsection is reserved for the (de)coupling of the dynamical systems equipped with irreversible characters, considered in the previous subsection. 

To this end, we first recall bicocycle double cross product Lie algebroid $\mathcal{A} {\hspace{.1cm}}_{\zeta\hspace{-.1cm}} \bowtie _\psi \mathcal{B} $ of Subsection \ref{sec-bicocycle}, with coordinates $(x^i,y^\alpha,y^a)$. In accordance with this, we shall assume the coordinates $(x^i,y_\alpha,y_a)$ on the dual bundle $\mathcal{A}^{\ast} {\hspace{.1cm}}_{\zeta\hspace{-.1cm}} \bowtie _\psi \mathcal{B}^{\ast}$, on which the Poisson structure was given in Subsection \ref{sec-dyn-bicocycle}. The bicocycle double cross product Poisson bracket, then, was given in \eqref{UP-Algoid-Poisson} in a coordinate-free way, and in \eqref{local Poisson on dual A} using coordinates. 

We then extend the bicocycle double cross product Lie algebroid with the real line to arrive at the extended product space $\mathcal{A}  {\hspace{.1cm}}_{\zeta\hspace{-.1cm}} \bowtie _\psi \mathcal{B} \times \mathbb{R}$ with local coordinates $(x^i,y^\alpha,y^a,z)$. The dual of this bundle is  $\mathcal{A}^{\ast} {\hspace{.1cm}}_{\zeta\hspace{-.1cm}} \bowtie _\psi \mathcal{B}^{\ast} \times \mathbb{R}$ with local coordinates $(x^i,y_\alpha,y_a,z)$.

\subsubsection*{Dissipative Hamiltonian Dynamics on the Dual BDCP Lie Algebroids.}

We extend the Poisson brackets in \eqref{local Poisson on dual A} following the approach defined in Subsection \ref{sec-irr-Lie}. More precisely, we merge \eqref{local Poisson on dual A} with the Jacobi brackets in \eqref{Pois-bra-A*xR}, where the extension with the action is examined. This itinerary map leads us to the following set of Jacobi brackets 
\begin{equation} \label{local Poisson on dual AxR+++}
\begin{split}
\{ x^{i}, x^{j} \}_{\mathcal{A}^{\ast} _\zeta \bowtie_\psi \mathcal{B}^{\ast}\times \mathbb{R} }  & = 0,\\
		  \{ z,x^{i}\}_{\mathcal{A}^{\ast} _\zeta \bowtie_\psi \mathcal{B}^{\ast}\times \mathbb{R} }&=0, \\ 
		\{z,y_\alpha\}_{\mathcal{A}^{\ast} _\zeta \bowtie_\psi \mathcal{B}^{\ast}\times \mathbb{R} }& =y_\alpha, \\ 
		 \{z,y_a\}_{\mathcal{A}^{\ast} _\zeta \bowtie_\psi \mathcal{B}^{\ast}\times \mathbb{R} }&=y_a,\\
\{x^{i}, y_{\alpha} \}_{\mathcal{A}^{\ast} _\zeta \bowtie_\psi \mathcal{B}^{\ast}\times \mathbb{R}} &  = (\mathfrak{a}_{\mathcal{A}})_{\alpha}^{i}, \\ 
\{x^{i}, y_{a} \}_{\mathcal{A}^{\ast} _\zeta \bowtie_\psi \mathcal{B}^{\ast}\times \mathbb{R}} &  = (\mathfrak{a}_{\mathcal{B}})_{a}^{i},\\ 
\{ y_{\alpha}, y_{\beta} \}_{\mathcal{A}^{\ast} _\zeta \bowtie_\psi \mathcal{B}^{\ast}\times \mathbb{R}}  &= \mathfrak{H}_{\alpha \beta}^{\gamma} y_{\gamma} + \mathfrak{Z}_{\alpha \beta}^{a} y_{a}, 
\\ 
\{ y_{a}, y_{b} \}_{\mathcal{A}^{\ast} _\zeta \bowtie_\psi \mathcal{B}^{\ast}\times \mathbb{R}}  &= \mathfrak{P}_{a b}^{\gamma} y_{\gamma} + \mathfrak{T}_{a b}^{d} y_{d},
\\ 
\{ y_{a}, y_{\beta} \}_{\mathcal{A}^{\ast} _\zeta \bowtie_\psi \mathcal{B}^{\ast}\times \mathbb{R}} &= \mathfrak{R}_{a \beta}^{\gamma} y_{\gamma} + \mathfrak{S}_{a \beta}^{d} y_{d}
\end{split}
\end{equation}
which we call the \emph{dissipative bicocycle double cross product Jacobi structure}. 

Accordingly, for a Hamiltonian function $\mathcal{H}=\mathcal{H}(x^i,y_\alpha,y_a,z)$ defined on the extended bundle $\mathcal{A}^{\ast} {\hspace{.1cm}}_{\zeta\hspace{-.1cm}} \bowtie _\psi \mathcal{B}^{\ast}\times \mathbb{R}$, the Hamilton's equations which we call the \emph{dissipative bicocycle double cross product Hamilton's equations}, is the following set of equations: 
\begin{equation}\label{dHamAlgebroid}
\begin{split}
	&\frac{dx^i}{dt}=(\mathfrak{a}_{\mathcal{A}})_{\alpha}^{i}\frac{\partial \mathcal{H}}{\partial y_{\alpha}} +(\mathfrak{a}_{\mathcal{A}})_{a}^{i}\frac{\partial \mathcal{H}}{\partial y_{a}},\\ & \frac{dy_\alpha}{dt}=-\Big( (\mathfrak{a}_{\mathcal{A}})_{\alpha}^{i}\frac{\partial \mathcal{H}}{\partial x^{i}}+ 
 \mathfrak{Z}_{\alpha \beta}^{a} y_{a}\frac{\partial \mathcal{H}}{\partial y_{\beta}}
 +\mathfrak{H}_{\alpha \beta}^{\gamma} y_{\gamma}\frac{\partial \mathcal{H}}{\partial y_{\beta}}
 - \mathfrak{R}_{b \alpha}^{\gamma} y_{\gamma}\frac{\partial \mathcal{H}}{\partial y_{b}}
 - \mathfrak{S}_{b \alpha}^{d} y_{d}\frac{\partial \mathcal{H}}{\partial y_{b}}
 + y_{\alpha}\frac{\partial \mathcal{H}}{\partial z}\Big) ,\\ & \frac{dy_a}{dt}=-\Big( (\mathfrak{a}_{\mathcal{B}})_{a}^{i}\frac{\partial \mathcal{H}}{\partial x^{i}}
 +
 \mathfrak{R}_{a \beta}^{\gamma} y_{\gamma}\frac{\partial \mathcal{H}}{\partial y_{\beta}}
 + 
  \mathfrak{S}_{a \beta}^{b} y_{b}\frac{\partial \mathcal{H}}{\partial y_{\beta}}
 + 
 \mathfrak{P}_{a b}^{\gamma} y_{\gamma}\frac{\partial \mathcal{H}}{\partial y_{b}}
 + \mathfrak{T}_{a b}^{d} y_{d}\frac{\partial \mathcal{H}}{\partial y_{b}}
 + y_{a}\frac{\partial \mathcal{H}}{\partial z}\Big) ,\\ & \frac{dz}{dt}=  y_{a}\frac{\partial \mathcal{H}}{\partial y_{a}} +y_{\alpha}\frac{\partial \mathcal{H}}{\partial y_{\alpha}}-\mathcal{H}    .
 \end{split}
\end{equation}

\subsubsection*{BDCP Lie-Poisson-Herglotz Equations.}

Let us reconsider the dissipative dynamics on Lie algebra setting, by assuming the base manifold to be a point, in which case the anchor map turns out to be trivial. Accordingly, the extended dual of the bicocycle double cross product space $\mathfrak{g} {\hspace{.1cm}}_{\zeta\hspace{-.1cm}} \bowtie _\psi \mathfrak{h}$ becomes $\mathfrak{g}^{\ast} {\hspace{.1cm}}_{\zeta\hspace{-.1cm}} \bowtie _\psi \mathfrak{h}^{\ast}\times \mathbb{R}$. Endowing the dual space with the coordinates $(y_\alpha,y_a,z)$, the Jacobi bracket in \eqref{local Poisson on dual AxR+++} takes the particular form
\begin{equation} \label{local Poisson on dual AxR++++}
\begin{split} 
		\{z,y_\alpha\}_{\mathfrak{g}^{\ast} _\zeta \bowtie_\psi \mathfrak{h}^{\ast}\times \mathbb{R} }& =y_\alpha, \\ 
		 \{z,y_a\}_{\mathfrak{g}^{\ast} _\zeta \bowtie_\psi \mathfrak{h}^{\ast}\times \mathbb{R} }&=y_a,\\
\{ y_{\alpha}, y_{\beta} \}_{\mathfrak{g}^{\ast} _\zeta \bowtie_\psi \mathfrak{h}^{\ast}\times \mathbb{R}}  &= \mathfrak{H}_{\alpha \beta}^{\gamma} y_{\gamma} + \mathfrak{Z}_{\alpha \beta}^{a} y_{a}, 
\\ 
\{ y_{a}, y_{b} \}_{\mathfrak{g}^{\ast} _\zeta \bowtie_\psi \mathfrak{h}^{\ast}\times \mathbb{R}}  &= \mathfrak{P}_{a b}^{\gamma} y_{\gamma} + \mathfrak{T}_{a b}^{d} y_{d},
\\ 
\{ y_{a}, y_{\beta} \}_{\mathfrak{g}^{\ast} _\zeta \bowtie_\psi \mathfrak{h}^{\ast}\times \mathbb{R}} &= \mathfrak{R}_{a \beta}^{\gamma} y_{\gamma} + \mathfrak{S}_{a \beta}^{d} y_{d}.
\end{split}
\end{equation}
Consequently, given a Hamiltonian function $h=h(y_\alpha,y_a,z)$ defined on  $\mathfrak{g}^{\ast}{\hspace{.1cm}}_{\zeta\hspace{-.1cm}} \bowtie _\psi \mathfrak{h}^{\ast}\times \mathbb{R}$, the \emph{bicocycle double cross product Lie-Poisson-Herglotz equations} are given by 
\begin{equation}\label{dHamAlgebra}
\begin{split}
& \frac{dy_\alpha}{dt}=-\Big(  \mathfrak{Z}_{\alpha \beta}^{a} y_{a}\frac{\partial h}{\partial y_{\beta}}+\mathfrak{H}_{\alpha \beta}^{\gamma} y_{\gamma}\frac{\partial h}{\partial y_{\beta}}
 - \mathfrak{R}_{b \alpha}^{\gamma} y_{\gamma}\frac{\partial h}{\partial y_{b}}
 - \mathfrak{S}_{b \alpha}^{d} y_{d}\frac{\partial h}{\partial y_{b}}
 + y_{\alpha}\frac{\partial h}{\partial z}\Big) ,\\ & \frac{dy_a}{dt}=-\Big(  \mathfrak{R}_{a \beta}^{\gamma} y_{\gamma}\frac{\partial h}{\partial y_{\beta}}+
  \mathfrak{S}_{a \beta}^{b} y_{b}\frac{\partial h}{\partial y_{\beta}} 
 + \mathfrak{P}_{a b}^{\gamma} y_{\gamma}\frac{\partial h}{\partial y_{b}}
 + \mathfrak{T}_{a b}^{d} y_{d}\frac{\partial h}{\partial y_{b}}
 + y_{a}\frac{\partial h}{\partial z}\Big) ,\\ & \frac{dz}{dt}=  y_{a}\frac{\partial h}{\partial y_{a}} +y_{\alpha}\frac{\partial h}{\partial y_{\alpha}}-h    .
 \end{split}
\end{equation}

\subsubsection*{Lagrangian  Dynamics on BDCP Lie Algebroids.} 

We shall now study the irreversible Lagrangian dynamics on the extended space $\mathcal{A}{\hspace{.1cm}}_{\zeta\hspace{-.1cm}} \bowtie _\psi \mathcal{B} \times \mathbb{R}$, equipped with the local coordinates $(x^i,y^\alpha,y^a,z)$ where $z$ denotes the coordinate on $\mathbb{R}$.  Accordingly, given a Lagrangian function $L=L(x^i,y^\alpha,y^a,z)$, we have two sets of equations as we have in reversible Lagrangian motion. The first set of equations is the one that determines the dynamics of the base space given by
\begin{equation}\label{LagLagLag}
\frac{dx^i}{dt} =  (\mathfrak{a}_{\mathcal{A}})^i_\alpha \,y^\alpha+(\mathfrak{a}_{\mathcal{B}})^i_a y^a.
\end{equation}
As it can be easily observed, this equation is the same with the reversible case. The second set of equations 
\begin{equation} \label{dLagAlgebroid}
\begin{split}
&  \frac{d}{dt}\frac{\partial L}{\partial y^\beta} =(\mathfrak{a}_{\mathcal{A}})^i_\beta\frac{\partial L}{\partial x^i} + \mathfrak{H}^\alpha_{\beta\gamma}y^\gamma\frac{\partial L}{\partial y^\alpha} -
\underbrace{\mathfrak{R}^\alpha_{d \beta }y^d \frac{\partial L}{\partial y^\alpha}}_{\text{\textit{rep.}}  ~ of ~ \mathcal{B} ~  on ~  \mathcal{A}} -
\underbrace{\mathfrak{S}^a_{d\beta }y^d\frac{\partial L}{\partial y^a}}_{\text{right \textit{rep.} of }\mathcal{A} \text{ on } \mathcal{B}}  \\
& \qquad \qquad\qquad\qquad +
\underbrace{\mathfrak{Z}^a_{\beta \gamma}y^\gamma\frac{\partial L}{\partial y^a}}_{\text{\textit{cocycle} term } \zeta } + \underbrace{\frac{\partial L}{\partial y^\beta}  \frac{\partial L}{\partial z}}_{\text{dissipative term}  },\\
&  \frac{d}{dt}\frac{\partial L}{\partial y^b} = (a_\mathcal{B})^i_b \frac{\partial L}{\partial x^i} +\mathfrak{T}^a_{bd}y^d\frac{\partial L}{\partial y^a} + 
\underbrace{\mathfrak{S}^a_{b\gamma}y^\gamma\frac{\partial L}{\partial y^a}}_{\text{\textit{rep.} of }\mathcal{A} \text{ on } \mathcal{B}}+ \underbrace{\mathfrak{R}^\alpha_{b\gamma}y^\gamma\frac{\partial L}{\partial y^\alpha}}_{\text{\textit{rep.}}  ~ of ~ \mathcal{B} ~  on ~  \mathcal{A}}  \\
& \qquad \qquad\qquad\qquad + \underbrace{\mathfrak{P}^\alpha_{bd}y^d\frac{\partial L}{\partial y^\alpha}}_{\text{\textit{cocycle} term } \psi } + \underbrace{\frac{\partial L}{\partial y^b}  \frac{\partial L}{\partial z}}_{\text{dissipative term}  }
\end{split}
\end{equation}
governs the dynamics on the fiber.

\subsubsection*{Lie Algebra Framework.} 

As a particular instance of the irreversible collective motion on the Lie algebroid setting, we shall now record the irreversible motion on the Lie algebra level. To this end, we consider the extension $\mathfrak{g}{\hspace{.1cm}}_{\zeta\hspace{-.1cm}} \bowtie _\psi \mathfrak{h}\times \mathbb{R}$ of the bicocycle double cross product Lie algebra $\mathfrak{g}{\hspace{.1cm}}_{\zeta\hspace{-.1cm}} \bowtie _\psi \mathfrak{h}$. The irreversible Lagrangian dynamics obtained on the extended space is now given by the bicocycle double cross product Euler-Poincar\'{e}-Herglotz equations.  Given a Lagrangian function $l=l(y^\alpha,y^a,z)$ on $\mathfrak{g}{\hspace{.1cm}}_{\zeta\hspace{-.1cm}} \bowtie _\psi \mathfrak{h}\times \mathbb{R}$, we have 
\begin{equation} \label{dLagAlgebra}
\begin{split}
&  \frac{d}{dt}\frac{\partial l}{\partial y^\beta} = \mathfrak{H}^\alpha_{\beta\gamma}y^\gamma\frac{\partial l}{\partial y^\alpha} -
\underbrace{\mathfrak{R}^\alpha_{d \beta }y^d \frac{\partial l}{\partial y^\alpha}}_{\text{\textit{rep.}}  ~ of ~ \mathfrak{h}  ~  on ~  \mathfrak{g} } -
\underbrace{\mathfrak{S}^a_{d\beta }y^d\frac{\partial l}{\partial y^a}}_{\text{\textit{rep.} of }\mathfrak{g}  \text{ on } \mathfrak{h} }   +
\underbrace{\mathfrak{Z}^a_{\beta \gamma}y^\gamma\frac{\partial l}{\partial y^a}}_{\text{\textit{cocycle} term } \zeta }\\
& \qquad \qquad\qquad\qquad  + \underbrace{\frac{\partial l}{\partial y^\beta}  \frac{\partial l}{\partial z}}_{\text{dissipative term}  },
 \\
&  \frac{d}{dt}\frac{\partial l}{\partial y^b} =  \mathfrak{T}^a_{bd}y^d\frac{\partial l}{\partial y^a} + 
\underbrace{\mathfrak{S}^a_{b\gamma}y^\gamma\frac{\partial l}{\partial y^a}}_{\text{\textit{rep.} of }\mathfrak{g}  \text{ on } \mathfrak{h} }+ \underbrace{\mathfrak{R}^\alpha_{b\gamma}y^\gamma\frac{\partial l}{\partial y^\alpha}}_{\text{\textit{rep.}}  ~ of ~ \mathfrak{h}  ~  on ~  \mathfrak{g} }  + \underbrace{\mathfrak{P}^\alpha_{bd}y^d\frac{\partial l}{\partial y^\alpha}}_{\text{\textit{cocycle} term } \psi } 
\\
& \qquad \qquad\qquad\qquad + \underbrace{\frac{\partial l}{\partial y^b}  \frac{\partial l}{\partial z}}_{\text{dissipative term}  }.
\end{split}
\end{equation}

\section{Conclusion}

The manuscript provides a comprehensive account of the (de)coupling phenomenon of Lie algebroids, and the dynamics on this geometry in the most general way through the bicocycle double cross product construction. This framework subsumes unified products, double cross products (matched pairs), cocycle extensions, and semi-direct product theories. The (de)coupling problem is addressed at the pure algebraic level in Section \ref{sec-unfLiealg}, along with an analysis from a dynamical point of view in Section \ref{dyndyn}. 

More precisely, in Theorem \ref{thm-muratcan} we have presented the compatibility conditions to couple two vector bundles under mutual interactions via twisted cocycles. Then in Theorem \ref{mezgi+} we showed that the bicocycle double cross product construction provides  the most tangible decomposition. Accordingly, the dynamics (dissipative or not) encoded by a vector bundle can be decomposed in simpler units for a better analysis. 

The bicocycle double cross product Hamiltonian dynamics in the Lie algebroid framework is written in \eqref{BDCPHamEqLoc}, while the bicocycle double cross product Lie-Poisson dynamics is in \eqref{BDCPHamEqLoc++}. The Lagrangian dynamics for bicocycle double cross products is given in \eqref{BDCP-EL1} and \eqref{BDCP-EL},  whereas the bicocycle double cross product Euler-Poincar\'{e} equation is given by \eqref{BDCP-EP}. The Herglotz type generalizations of all these dynamical models are listed in Subsection \ref{noveldyn}. 

As a future work, we wish to address the (de)coupling problem on Filippov $3$-Lie algebroids, in an attempt to study collective Nambu motion.


\bibliographystyle{abbrv}
\bibliography{references}{}

\end{document}